\documentclass{amsart}

\newtheorem{theorem}{Theorem}[section]
\newtheorem{proposition}[theorem]{Proposition}
\newtheorem{corollary}[theorem]{Corollary}
\newtheorem{lemma}[theorem]{Lemma}

\newtheorem{preremark}[theorem]{Remark}
\newtheorem{predefinition}[theorem]{Definition}
\newtheorem{preexample}[theorem]{Example}
\newtheorem{prenotation}[theorem]{Notation}
\newtheorem{preconjecture}[theorem]{Conjecture}

\newenvironment{remark}{\begin{preremark}\rm}{\end{preremark}}
\newenvironment{definition}{\begin{predefinition}\rm}{\end{predefinition}}
\newenvironment{example}{\begin{preexample}\rm}{\end{preexample}}
\newenvironment{notation}{\begin{prenotation}\rm}{\end{prenotation}}

\def\LT{\mathop{\rm LT}\nolimits}

\def\Mat{\mathop{\rm Mat}\nolimits}

\def\Supp{\mathop{\rm Supp}\nolimits}
\def\Spec{\mathop{\rm Spec}\nolimits}
\def\Hilb{\mathop{\rm Hilb}\nolimits}

\def\BO{\mathbb{B}_{\mathcal{O}}}
\def\BP{\mathbb{B}_{\mathcal{P}}}
\def\CO{\mathbb{C}_{\mathcal{O}}}

\def\yi{\mathbf{y}^{(i)}}
\def\y1{\mathbf{y}^{(1)}}
\def\yd{\mathbf{y}^{{}^{\hbox{\bf .}}}}
\def\ymu{\mathbf{y}^{(\mu)}}
\def\mod{\mathop{\rm mod}\,}

\let\epsilon=\varepsilon
\let\rho=\varrho
\let\phi=\varphi
\let\To=\longrightarrow

\def\tfrac #1#2{{\textstyle\frac{#1}{#2}}}

\def\cocoa{\mbox{\rm
 C\kern-.13em o\kern-.07 em C\kern-.13em o\kern-.15em A}}
\def\apcocoa{\mbox{\rm
A\kern-0.13em p\kern -0.07em C\kern-.13em o\kern-.07 em C\kern-.13em
o\kern-.15em A}}

\begin{document}

\title{The Geometry of Border Bases}

%    Information for first author
\author{Martin Kreuzer}
\address{Fakult\"at f\"ur Informatik und Mathematik, Universit\"at
Passau,
D-94030 Passau, Germany}
\email{Martin.Kreuzer@uni-passau.de}

%    Information for second author
\author{Lorenzo Robbiano}
\address{Dipartimento di Matematica, Universit\`a di Genova, Via
Dodecaneso 35,
I-16146 Genova, Italy}
\email{robbiano@dima.unige.it}

\date{\today}

\begin{abstract}
In this paper we continue the study of the border basis scheme
we started in~\cite{KR3}. The main topic is the construction
of various explicit flat families of border bases. To begin with, 
we cover the punctual Hilbert scheme $\Hilb^\mu(\mathbb{A}^n)$ 
by border basis schemes and work out the base changes. 
This enables us to control
flat families obtained by linear changes of coordinates. Next we
provide an explicit construction of the principal component of the
border basis scheme, and we use it to find flat families of maximal
dimension at each radical point. Finally, we connect radical points to each
other and to the monomial point via explicit flat families on the
principal component.
\end{abstract}


\maketitle

%%%%%%%%%%%%%%%%%%%%%%%%%%%%%%%%%%%%%%
%  Section 1: Introduction
%%%%%%%%%%%%%%%%%%%%%%%%%%%%%%%%%%%%%%

\section{Introduction}

\begin{flushright}
{\it Nobody goes there anymore}\\
{\it because it's too crowded.}\\
{\rm (Yogi Berra)}\\
\end{flushright}

Border bases schemes have recently become an active area
of research, as evidenced by a list of references which is 
getting quite crowded, e.g.~\cite{ABM}, \cite{GLS}, \cite{Hu1},
\cite{Hu4}, \cite{KPR}, \cite{KR3}, \cite{Led}, \cite{Lee}, 
and~\cite{R}, to mention just a few contributions of the
last years. What is the reason for this spurt of activity?

In our opinion there are several reasons. One of them is that border
bases enjoy a degree of numerical stability which, in contrast,
Gr\"obner bases don't. This has proven useful for dealing with
empirical polynomials constructed from measured data (see for
instance~\cite{KPR} and~\cite{S}) and has even led to actual
industrial applications. But the most relevant aspect for our
topic is that border basis schemes provide a very concrete
and easily accessible way to parametrize 0-dimensional polynomial
ideals. They can be viewed as open affine subschemes of the corresponding
Hilbert schemes which can be described by simple, explicit polynomial
equations (see for instance~\cite{KR3} and~\cite{MS}).

This brings us to our first contribution: by constructing
explicit matrices describing the change of basis between
one border basis scheme and another, we obtain a direct 
construction of the punctual Hilbert scheme $\Hilb^\mu(\mathbb{A}^n)$
(see~\cite{HS} and~\cite{Hu2}) which
uses neither A.\ Grothendieck's Grassmannian variety technique
nor any arguments involving representation of functors.

In the paper~\cite{R} it was shown that in some cases border
bases schemes can be described via suitable Gr\"obner basis
schemes. This is an important fact, since the description of
Gr\"obner basis schemes requires fewer indeterminates
and fewer equations, and motivates the strategy
used in Section~1 to treat the cases of border basis
and Gr\"obner basis schemes simultaneously. 

Another {\it driving force} for writing this paper was the search for
suitable flat deformations of border bases.
We have seen in~\cite{KR3} that flat deformations of border bases
are the same as rational curves on the border basis scheme~$\BO$.
Therefore we want to construct explicit rational curves on~$\BO$.
However, in contrast to the Hilbert scheme, we have the problem
that a flat family whose general fibre is an ideal corresponding to
a point in~$\BO$, i.e.\ an ideal having an~$\mathcal{O}$-border
basis, can have a special fibre which doesn't (see
for instance~\cite{KR3}, Example~3.9). Constructing many flat
families of border bases could be one way to attack the
hitherto unsolved problem of the connectedness of the
border basis scheme (see~\cite{R}, Question~2). 
In various parts of this paper we construct
a number of such flat families, in addition to the ones
we found by homogenization in~\cite{KR3}:
families obtained from linear changes of coordinates, 
from local parametrizations of the principal component, and
from distractions.

\smallskip
Let us describe the contents of the paper in more detail.
In Section 2 we introduce pseudo order ideals,
pseudo borders and pseudo border bases 
(see Definition~\ref{pseudodef}).
Using them, we are not only able to treat the cases of
order ideals with their borders and of $\sigma$-cornercuts with their
corners simultaneously, but also their isomorphic images under 
a linear change of coordinates.
After constructing a moduli space for pseudo border bases
(cf.\ Prop.~\ref{OBorderScheme}) and showing that it
comes equipped with a universal flat family and has the 
expected properties (cf.\ Prop.~\ref{PseudoUF}), we
arrive at our first main result.
In Theorem~\ref{basechange} we construct an explicit
isomorphism between the open subsets of two pseudo
border basis schemes corresponding to the ideals
which have pseudo border bases with respect to both
pseudo order ideals.

As mentioned above, this yields an explicit description 
of the glueing of border basis schemes necessary
to build the corresponding punctual Hilbert scheme (see
Remark~\ref{constrHilb}). Another application of the theorem 
is the possibility to characterize when a linear change of
coordinates produces an ideal which has again a pseudo
border basis with respect to the same pseudo order ideal
(see Prop.~\ref{lc_represent}). Thus we can use linear changes 
of coordinates to construct explicit flat families of border bases
(see Proposition~\ref{linchange}). Moreover, we show that, in the 
$\mathcal{O}$-border basis setting, generic linear changes of 
coordinates lead to ideals which have $\mathcal{O}$-border bases again
(see Corollary~\ref{no-gin-cor}) and conclude quippingly
that in border basis theory there in {\it no gin}.

In Section~3 we start our exploration of the principal component 
of the border basis scheme (see Definition~\ref{defprincomp}). 
The first step is Theorem~\ref{princicomp} where we provide 
explicit equations defining this scheme.
Next we show in Proposition~\ref{COgens} that our construction
yields the same result as the one given in~\cite{Ha1} and~\cite{Ha2},
but uses a much smaller number of algebra generators.
This has the obvious advantage that one can turn our description
into an algorithm for computing the vanishing ideal of the principal
component (see Prop.~\ref{PCcomp}), and that one can use this algorithm
to check whether a given border basis scheme is irreducible.

In Section~4 we use the principal component~$\CO$ of~$\BO$ to
construct more explicit flat families of
border bases. More precisely, we construct explicit local 
parameters at a radical point of~$\CO$. The idea of this
construction is to use the complete intersection representation
of a radical ideal provided by the Shape Lemma (c.f.~\cite{KR2},
Theorem~3.7.25) and to apply the techniques of Section~1 to it.
The resulting Theorem~\ref{localparams} not only recovers
well-known facts, e.g.\ that~$\CO$ is a rational variety and
non-singular at its radical points,
but it gives us an explicit parametrization of an open neighborhood
of every radical point. We use this theorem in several ways:
to connect two radical points on~$\CO$ via a sequence of two explicit
flat families (cf.\ Remark~\ref{connect2rad}), and to combine these 
families with a distraction to connect every radical point of~$\CO$ 
to the monomial point (cf.\ Remark~\ref{connect2mon}).

\smallskip
Unless explicitly stated otherwise, we use the notation and the
definitions introduced in~\cite{KR1} and~\cite{KR2}.

\bigskip
%%%%%%%%%%%%%%%%%%%%%%%%%%%%%%%
%
%  Section 2: Change of Basis
%
%%%%%%%%%%%%%%%%%%%%%%%%%%%%%%%

\section{Change of Basis}

Let $K$ be a field, let $P=K[x_1,\dots,x_n]$,
let $I\subset P$ be a 0-dimensional ideal, and let
$\mu=\dim_K(P/I)$.
As mentioned in the introduction, it is our goal to
construct explicit flat families of ideals
having~$I$ as one of their fibers.
A natural approach is to perform linear changes of
coordinates, i.e.\ $K$-algebra isomorphisms
$\phi:P\longrightarrow P$ mapping the indeterminates to
(not necessarily homogeneous) linear polynomials.
If the ideal~$I$ has a border basis with respect to
some order ideal, the same is not always true for the
ideal $\phi(I)$. Therefore one of the ideas of
the following construction is to keep track when
a linear change of coordinates preserves the 
property that~$I$ has a border basis with respect to
a given order ideal.

Recall that a finite set of terms $\mathcal{O}$ in~$\mathbb{T}^n$
is called an {\bf order ideal} if it is closed under forming 
divisors, i.e.\ if $t\in\mathcal{O}$ and $t'\mid t$ imply
$t'\in\mathcal{O}$. The set of terms $\partial\mathcal{O}=
(x_1\mathcal{O}\cup \cdots \cup x_n\mathcal{O})\setminus
\mathcal{O}$ is called the {\bf border} of~$\mathcal{O}$.
The definition of an order ideal implies that the set
$\mathbb{T}^n\setminus \mathcal{O}$ is a monomial ideal.
We denote the set of monomial generators of this monomial
ideal by~$c\mathcal{O}$ and call it the {\bf corner set}
of~$\mathcal{O}$. Let $\sigma$ be a term ordering. 
The order ideal~$\mathcal{O}$ is called
a {\bf $\sigma$-cornercut} if $b >_\sigma t$ for all $b\in c\mathcal{O}$ 
and all $t\in \mathcal{O}$. Notice that this implies $b >_\sigma t$
for all $b\in \partial\mathcal{O}$ and all $t\in \mathcal{O}$.

The following definition is manufactured in such
a way that we can treat the cases of the border basis scheme and
the Gr\"obner basis scheme simultaneously.

\begin{definition}\label{pseudodef}
Let $\phi: P\longrightarrow P$ be a linear change of coordinates,
and let~$\sigma$ be a term ordering.

\begin{enumerate}
\item[(a)] Let~$\mathcal{P}$ and~$b\mathcal{P}$ be sets of
polynomials  in~$P$. Then~$\mathcal{P}$ is called
a {\bf pseudo order ideal} and~$b\mathcal{P}$ is called
the {\bf pseudo border}  of~$\mathcal{P}$ if one of the following
two cases occurs:
\begin{enumerate}
\item[(i)] $\mathcal{P}$ is the image of an order
ideal~$\mathcal{O}$ under~$\phi$, and~$b\mathcal{P}$ is the
corresponding image of the border of~$\mathcal{O}$;

\item[(ii)] $\mathcal{P}$ is the image of a
$\sigma$-cornercut~$\mathcal{O}$ under~$\phi$,
and~$b\mathcal{P}$ is the corresponding image of the
corner set~$c\mathcal{O}$.
\end{enumerate}

\item[(b)] Let $\mathcal{P}$ be a pseudo order ideal,
and let~$I$ be an ideal in~$P$ such that the residue classes
of the elements of~$\mathcal{P}$ form
a $K$-vector space basis of~$P/I$. In this case we simply say
that~$\mathcal{P}$ is a {\bf basis modulo~$I$}.

\item[(c)] Let $\mathcal{P}=\{t_1,\dots,t_\mu\}$ be a 
pseudo order ideal in~$P$, let $b\mathcal{P}=\{b_1,\dots,b_\nu\}$ 
be its pseudo border, and for $j=1,\dots,\nu$ let
$g_j=b_j - \sum_{i=1}^\mu \gamma_{ij}t_i$ with $\gamma_{ij}\in K$.
Then the set $G=\{g_1,\dots,g_\nu\}$
is called a {\bf pseudo $\mathcal{P}$-border prebasis}.

\item[d)] A pseudo $\mathcal{P}$-border prebasis
$G=\{g_1,\dots,g_\nu\}$ is called
a {\bf pseudo $\mathcal{P}$-border basis} if~$\mathcal{P}$
is a basis modulo the ideal $(g_1,\dots,g_\nu)$.

\item[(e)] Let $\alpha=\#(b\mathcal{P})$, and let
$C=(c_{ij})$ be a matrix of indeterminates of size $\mu\times\alpha$.
For $j=1,\dots,\alpha$, we form the polynomials
$g_j = b_j-\sum_{i=1}^\mu  c_{ij}t_i$. Then $G=\{g_1,\dots,g_\alpha\}$
is called the {\bf generic pseudo $\mathcal{P}$-border prebasis}.

\end{enumerate}
\end{definition}

In the following we shall assume the setting and notation
of this definition. We start our investigation 
with the following fact.

\begin{proposition}\label{OBorderScheme}
Let $\mathcal{P}=\{t_1,\dots,t_\mu\}$ be a pseudo order
ideal in~$P$, let $b\mathcal{P}=\{b_1,\dots,b_\alpha\}$
be its pseudo border, and let~$C=(c_{ij})$ be a
matrix of indeterminates of size $\mu\times\alpha$.
There exists an ideal~$I(\BP)$ in $K[c_{11},\dots, c_{\mu\alpha}]$
such that the ring $B_\mathcal{P}=K[c_{11}, \dots, c_{\mu\alpha}]/I(\BP)$
is the coordinate ring of an affine scheme $\BP$ 
whose $K$-rational points are in one-to-one correspondence
with the ideals~$I$ in~$P$ for which~$\mathcal{P}$
is a basis modulo~$I$.
\end{proposition}

\begin{proof}
If $\mathcal{P}=\mathcal{O}$ is an order ideal of terms
and~$b\mathcal{P} = \partial\mathcal{O}$
its border, the claim follows from~\cite{KR2}, Theorem 6.4.30.
If $\mathcal{P}=\mathcal{O}$ is a $\sigma$-cornercut for some 
term ordering~$\sigma$ and $b\mathcal{P} = c\mathcal{O}$, the claim
follows from~\cite{R}, Proposition 3.11.
Given an ideal $I\subset P$ such
that~$\mathcal{O}$ is a basis modulo~$I$, 
let~$C_I$ be the matrix obtained by substituting the
entries~$c_{ij}$ of~$C$ with the coordinates of the
point in the scheme~$\BO$ corresponding to~$I$.
We observe that, in both cases, we have
$$
b\mathcal{O} = \mathcal{O}\cdot C_I \quad  {\rm mod\ } I \eqno{(1)}
$$

Next, let~$\mathcal{O}$ be an order ideal satisfying one of the
preceding two conditions, let $\phi: P\longrightarrow P$ be a linear 
change of coordinates, and let $\mathcal{P} = \varphi(\mathcal{O})$.
By definition, we have $b\mathcal{P}= \varphi(b\mathcal{O})$. We
apply~(1) to~$\mathcal{O}$ and $b\mathcal{O}$, and we get
$$
\phi(b\mathcal{O}) = \phi(\mathcal{O})\cdot C_I \quad {\rm mod\ } \phi(I)
\eqno (2)
$$
for all ideals~$I$ in~$P$ such that~$\mathcal{O}$ is a basis modulo~$I$,
and therefore
$$
b\mathcal{P} = \mathcal{P}\cdot C_I \quad {\rm mod\ } \phi(I) \eqno (3)
$$

Now let~$J$ be an ideal in~$P$ such that~$\mathcal{P}$ is a basis modulo~$J$.
Then $\phi^{-1}(J)$ is an ideal such that~$\mathcal{O}$ is a basis
modulo $\phi^{-1}(J)$, and we can use~(2) and~(3).
The fact that $J = \phi(\phi^{-1}(J))$ implies
$b\mathcal{P}=\mathcal{P}\cdot C_{\phi^{-1}(J)}$ modulo~$J$.
Hence, if we define $\BP=\BO$ and $I(\BP)=I(\mathbb{B}_{\mathcal{O}})$,
the ideals~$J$ in~$P$ such that~$\mathcal{P}$ is a basis modulo~$J$
correspond one-to-one to the ideals $I=\phi^{-1}(J)$
in~$P$ such that~$\mathcal{O}$ is a basis modulo~$I$.
\end{proof}

In the setting of this proposition, we introduce
further terminology.

\begin{definition}
As above, let $\mathcal{P}$ be a pseudo order ideal,
$b\mathcal{P}$ its pseudo border, and $\alpha=\#(b\mathcal{P})$.

\begin{enumerate}
\item[(a)] The scheme $\BP$ is called the
{\bf $\mathcal{P}$-basis scheme}.

\item[(b)] Given an ideal~$I$ in the polynomial ring~$P$ 
such that~$\mathcal{P}$ is
a basis modulo~$I$, we write $b\mathcal{P}=\mathcal{P}\cdot
C_I$ modulo~$I$. Then the matrix $C_I\in\Mat_{\alpha\mu}(K)$ and
the point $c_I$ whose coordinates are the entries of $C_I$ are said to
{\bf represent~$I$} in~$\mathbb{B}_\mathcal{P}$.

\item[(c)] Let $G=\{g_1,\dots,g_\alpha\}$ be the generic pseudo
$\mathcal{P}$-border prebasis, and let
$$
U_\mathcal{P} \;=\; K[x_1, \dots, x_n, c_{11},\dots, c_{\mu\alpha}]/
\big( I(\BP) + (g_1,\dots,g_\alpha)\big)
$$
Then the natural homomorphism of $K$-algebras
$\Phi: B_\mathcal{P} \To U_\mathcal{P}$
is called the {\bf universal $\mathcal{P}$-basis family}.
\end{enumerate}
\end{definition}

\begin{remark}
Let us point out one fact that follows from the proof of the preceding
proposition: given a linear change of coordinates
$\phi: P\longrightarrow P$, the matrix $C_I$ which 
represents an ideal~$I$  in~$\mathbb{B}_\mathcal{O}$ 
also represents~$\phi(I)$ in~$\mathbb{B}_{\phi(\mathcal{O})}$.
\end{remark}

As in the usual border basis theory, a central property
of the universal family is that it is free with basis~$\mathcal{P}$.
This is the main part of the following proposition.

\begin{proposition}\label{PseudoUF}
As above, let $\mathcal{P}$ be a pseudo order ideal, 
let~$b\mathcal{P}$ be its pseudo
border, let $G=\{g_1,\dots,g_\alpha\}$ be the
generic pseudo $\mathcal{P}$-border prebasis, and let
${\Phi: B_\mathcal{P} \To U_\mathcal{P}}$ be the universal
$\mathcal{P}$-basis family.

\begin{enumerate}
\item[(a)] The residue classes of the elements of $\mathcal{P}$ are
a $B_\mathcal{P}$-module basis of~$U_\mathcal{P}$.

\item[(b)] Via $\Phi$, the ring~$B_\mathcal{P}$ is a
free summand of~$U_\mathcal{P}$.
In particular, the map~$\Phi$ is injective and~$B_\mathcal{P}$
can be seen as a subring of~$U_\mathcal{P}$.

\item[(c)] The rewriting rules given by the tuple $(g_1,\dots,g_\alpha)$
yield an explicit division algorithm with the following
properties: for every polynomial~$f$ in
$K[x_1,\dots,x_n, c_{11},\dots,c_{\mu\alpha}]$, it
produces a polynomial~$f'$ which is a
linear combination of the elements in~$\mathcal{P}$
with coefficients in $K[c_{11},\dots, c_{\mu\alpha}]$.
The residue classes of these coefficients in~$B_\mathcal{P}$
are uniquely determined and do not depend on the
ordering of $(g_1,\dots,g_\nu)$.
\end{enumerate}
\end{proposition}

\begin{proof}
First we show~(a). Let $\mathcal{O}$ be an order ideal in~$\mathbb{T}^n$
and $\phi:P\longrightarrow P$ a linear change of coordinates such that
$\mathcal{P}=\phi(\mathcal{O})$. By~\cite{KR3}, Theorem 3.4 and~\cite{R},
Theorem 2.9, respectively, the set~$\mathcal{O}$ is
a $B_{\mathcal{O}}$-module basis of~$U_{\mathcal{O}}$.
Now the claim follows from the fact that we used $I(\BP)=
I(\BO)$ in Proposition~\ref{OBorderScheme}.

Since~(b) follows immediately from~(a), it remains to prove~(c).
We denote the extension of~$\phi$ to
$P[c_{ij}]$ by $\tilde\phi$, and
again we write $\mathcal{P}=\phi(\mathcal{O})$ with an order
ideal~$\mathcal{O}$. To define the algorithm for
dividing~$f$ by $(g_1,\dots,g_\alpha)$, we use the usual border
division algorithm (cf.~\cite{KR2}, Proposition 6.4.11) to divide
$\tilde \phi^{-1}(f)$ by $(\tilde \phi^{-1}(g_1),\dots,\tilde\phi^{-1}(g_\alpha))$.
At the end we apply~$\tilde\phi$ to the resulting representation
of~$\tilde\phi^{-1}(f)$ as a $K[c_{ij}]$-linear combination of
the elements of~$\mathcal{O}$ and obtain the desired
$K[c_{ij}]$-linear combination of elements of~$\mathcal{P}$.
The uniqueness of this representation is a consequence of~(a).
\end{proof}

Next we let~$I$ be an ideal in~$P$ such that two
pseudo order ideals~$\mathcal{P}$ and~$\mathcal{P'}$ are bases modulo~$I$.
Suppose that~$I$ is represented by a matrix~$C_I$ in $\BP$
and a matrix $D_I$ in $\mathbb{B}_\mathcal{P'}$.
What is the relation between $C_I$ and $D_I$?
In the following we examine this question.

\begin{remark}
Let $\mathcal{P}$ and $\mathcal{P}'$ be two pseudo order ideals
for which $\mu=\#\mathcal{P}=\#\mathcal{P}'$.
The scheme $\BP$ parametrizes ideals~$I$ in~$P$ such that
$\mathcal{P}$ is a basis modulo~$I$. It contains a Zariski-open
subset $\mathbb{B}_{\mathcal{P}, \mathcal{P}'}$ which parametrizes
the ideals~$I$ such that also~$\mathcal{P}'$ is a basis modulo~$I$.
Similarly, there is an open subset~$\mathbb{B}_{\mathcal{P}',\mathcal{P}}$
of~$\mathbb{B}_{\mathcal{P}'}$ which
parametrizes the ideals~$I$ for which both~$\mathcal{P}$ and~$\mathcal{P}'$
are bases modulo~$I$.

It is known that $\mathbb{B}_\mathcal{P}$
and~$\mathbb{B}_\mathcal{P}'$ are open subschemes of the same punctual
Hilbert scheme $\Hilb^\mu(\mathbb{A}^n)$
(see for instance~\cite{MS}, Chapter 18).
Their intersection includes a non-empty open subset of the
principal component (i.e.\ the component corresponding to radical ideals,
cf.\ Section~\ref{The Principal Component of the Border Basis Scheme}).
Consequently, the open subsets 
$\mathbb{B}_{\mathcal{P},\mathcal{P}'}$~of $\Hilb^\mu(\mathbb{A}^n)$, 
which are equal by definition, are not empty.
\end{remark}

In the following, we let $\mathcal{P}$ and~$\mathcal{P}'$ be two
pseudo order ideals such that $\mu=\#\mathcal{P}=\#\mathcal{P}'$.
Let $\alpha=\#(b\mathcal{P})$ and $\alpha'=\#(b\mathcal{P}')$
be the cardinalities of their pseudo borders, let
$C=(c_{ij})$ be a matrix of indeterminates of size $\mu\times\alpha$,
and let $D=(d_{ij})$ be a matrix of indeterminates of size
$\mu\times\alpha'$. According to Proposition~\ref{PseudoUF}.c,
there exist matrices~$M_C$ and~$N_C$ over~$K[c_{ij}]$
such that
$$
\mathcal{P}' = \mathcal{P} \cdot M_C  \qquad  \qquad
b\mathcal{P'} = \mathcal{P} \cdot N_C \eqno{(4)}
$$
are matrix equalities which hold over~$U_\mathcal{P}$.
Similarly, there exist matrices~$M_D$ and $N_D$
over~$K[d_{ij}]$ such that
$$
\mathcal{P}= \mathcal{P'} \cdot M_D  \qquad  \qquad
b\mathcal{P} = \mathcal{P'} \cdot N_D \eqno{(5)}
$$
are matrix equalities which hold over~$U_\mathcal{P'}$.

\begin{theorem}{\rm (Base Change for Pseudo Border Basis Schemes)}
\label{basechange}\\
Assume that we are in the setting described above.

\begin{enumerate}
\item[(a)] The set $\mathbb{B}_{\mathcal{P}, \mathcal{P}'}$
is the open subset of~$\BP$ defined by $\det (M_C)\ne 0$, and the
set $\mathbb{B}_{\mathcal{P}', \mathcal{P}}$ is the open subset
of~$\mathbb{B}_\mathcal{P}'$ defined by~$\det (M_D) \ne 0$.

\item[(b)] The natural maps defining the identity map between
$\mathbb{B}_{\mathcal{P}, \mathcal{P'}}$
and~$\mathbb{B}_{\mathcal{P'}, \mathcal{P}}$
in terms of their respective systems of coordinates are given 
parametrically by
$$
D = M_C^{-1}\cdot N_C   {\quad  \rm and\quad}
C = M_D^{-1}\cdot N_D
$$
\end{enumerate}
\end{theorem}

\begin{proof}
Claim~(a) follows immediately from~(4) and~(5).
Now we prove claim~(b). By the definition of the
generic pseudo border prebases, we have the equalities
$$
b\mathcal{P} = \mathcal{P} \cdot C  \qquad \hbox{\rm and}\qquad
b\mathcal{P'} = \mathcal{P'} \cdot D \eqno{(6)}
$$
In the following we work over the open set where both~$M_C$ and~$M_D$
are invertible, i.e.\ where both systems of coordinates $(c_{ij})$
and $(d_{ij})$ apply.
Combining the second equality in~$(6)$ with the first in~$(4)$, we get
$$
b\mathcal{P}' = \mathcal{P}\cdot M_C\cdot D
$$
Now we compare this to the second equality in~$(4)$ and use the
uniqueness implied by the fact that~$\mathcal{P}$
is a $B_{\mathcal{P}}$-basis of~$U_{\mathcal{P}}$ to get
$$
M_C\cdot D = N_C
$$
This implies the first claim in~(b). The second claim follows by
interchanging the roles of~$\mathcal{P}$ and~$\mathcal{P}'$.
\end{proof}

The formulas given in this theorem can be used as follows to 
give an explicit construction for the punctual Hilbert scheme. 

\begin{remark}(Glueing Border Basis Schemes)\label{constrHilb}\\
Given an integer $\mu>0$, it is well-known that there exists a
scheme, called the {\it punctual Hilbert scheme} and denoted
by $\Hilb^\mu(\mathbb{A}^n)$ which parametrizes all 0-dimensional
ideals~$I$ in~$P$ such that $\dim_K(P/I)=\mu$.
For an introduction to punctual Hilbert schemes, see for instance~\cite{MS}
and its bibliography. Here we just want to point out that
Theorem~\ref{basechange} allows to construct~$\Hilb^\mu(\mathbb{A}^n)$
very explicitly.

To explain the method, we use the following example.
Let $n=2$ and $\mu=4$, i.e.\ we want to parametrize ideals
which correspond to 0-dimensional subschemes of~$\mathbb{A}^2$
of length four. There exist exactly five order ideals which can
serve as a $K$-basis of~$P/I$, namely
\begin{eqnarray*}
\mathcal{O}_1 &=& \{1,\, x,\, x^2,\, x^3\}, 
\quad \mathcal{O}_2 = \{1,\, y,\, y^2,\, y^3\},
\quad  \mathcal{O}_3 = \{1,\, x,\, y,\, x^2\}, \\
\mathcal{O}_4 &=& \{1,\, x,\, y,\, y^2\}, \quad \hbox{\rm and}\quad
\mathcal{O}_5 = \{1,\, x,\, y,\, xy\}.
\end{eqnarray*}

Let $I_1 = (x^4,y)$, $I_2 =(x,y^4)$, $I_3 =(xy,y^2,x^3)$,
$I_4 = (x^2,xy,y^3)$, and $I_5 =(x^2,y^2)$. 
It is clear that, for $i=1,\dots,5$, among the five
order ideals the only basis modulo~$I_i$ is $\mathcal{O}_i$.
Therefore all the border basis schemes of all five order ideals 
are needed to cover $\Hilb^4(\mathbb{A}^2)$ with affine open sets.
The scheme $\Hilb^4(\mathbb{A}^2)$ can be constructed
explicitly by glueing the schemes $\mathbb B_{\mathcal{O}_1},\dots,
\mathbb{B}_{\mathcal{O}_5}$ via the isomorphisms given in
Theorem~\ref{basechange}.

Furthermore, we note that, for $i=1,\dots,4$, the set $\mathcal{O}_i$ 
is a $\sigma_i$-cornercut for a suitable term ordering~$\sigma_i$.
However, this is not the case for~$\mathcal{O}_5$. Consequently, we have 
${\mathbb B_{{\mathcal{O}_i}} = \mathbb G_{\sigma_i,\mathcal{O}_i}}$
for $i=1,\dots,4$ (see~\cite{R}, Proposition~3.11).
On the other hand, a direct calculation yields
$\dim(\mathbb B_{{\mathcal{O}_5}}) = 8$, as expected, but
$\dim(\mathbb G_{\sigma,\mathcal{O}_5})\le 7$
for every term ordering~$\sigma$. This implies that
we cannot cover $\Hilb^4(\mathbb{A}^2)$ with open sets associated to
Gr\"obner basis schemes. Using Gr\"obner basis schemes, we merely get a
stratification of the Hilbert scheme.
\end{remark}

The following example illustrates the theorem.

\begin{example}
In the ring $P=\mathbb{Q}[x,y]$ we consider the two order ideals
$\mathcal{O}=\{1,y,x,xy\}$ and $\mathcal{O}'=\{1,x,x^2,x^3\}$.
Our goal is to find the transformation matrices mentioned in the
preceding theorem and to verify the equations given in its proof.
(We have ordered all sets and tuples of terms according to
{\tt DegRevLex}.)

First of all, we represent~$\mathcal{O}'$ in terms of~$\mathcal{O}$,
modulo the generic $\mathcal{O}$-border basis $G=\{g_1,g_2,g_3,g_4\}$
where
\begin{eqnarray*}
g_1 &=& y^2 -c_{11} -c_{21}y -c_{31}x -c_{41}xy\\
g_2 &=& x^2 -c_{12} -c_{22}y -c_{32}x -c_{42}xy\\
g_3 &=& xy^2 -c_{13} -c_{23}y -c_{33}x -c_{43}xy\\
g_4 &=& x^2y -c_{14} -c_{24}y -c_{34}x -c_{44}xy
\end{eqnarray*}
The result is
$$
\mathcal{O}' \;=\; \mathcal{O} \cdot M_C \;=\; \mathcal{O}
\cdot \begin{pmatrix}
1 & 0 & c_{12} & c_{12}c_{32}+c_{14}c_{42}\\
0 & 0 & c_{22} & c_{22}c_{32}+c_{24}c_{42}\\
0 & 1 & c_{32} & c_{12}+c_{32}^2+c_{34}c_{42}\\
0 & 0 & c_{42} & c_{22}+c_{32}c_{42}+c_{42}c_{44}
\end{pmatrix}
$$
Similarly, we represent $\partial\mathcal{O}'$ in terms of $\mathcal{O}$
and find
$$
\partial\mathcal{O}' \;=\; \mathcal{O} \cdot N_C \;=\; \mathcal{O}
\cdot \begin{pmatrix}
0 & 0 & c_{14} &  c_{12}c_{34} + c_{14}c_{44} & h_1\\
1 & 0 & c_{24} &  c_{22}c_{34} + c_{24}c_{44} & h_2\\
0 & 0 & c_{34} &  c_{32}c_{34} + c_{34}c_{44} + c_{14} & h_3\\
0 & 1 & c_{44} &  c_{34}c_{42} + c_{44}^2 + c_{24} & h_4
\end{pmatrix}
$$
where
\begin{eqnarray*}
h_1 &=& c_{12}c_{32}^2 + c_{14}c_{32}c_{42} + c_{12}c_{34}c_{42} + c_{14}c_{42}c_{44}
  + c_{12}^2 + c_{14}c_{22}\\
h_2 &=& c_{22}c_{32}^2 + c_{24}c_{32}c_{42} + c_{22}c_{34}c_{42} + c_{24}c_{42}c_{44}
  + c_{12}c_{22} + c_{22}c_{24}\\
h_3 &=&  c_{32}^3 + 2c_{32}c_{34}c_{42} + c_{34}c_{42}c_{44} + 2c_{12}c_{32}
  + c_{22}c_{34} + c_{14}c_{42}\\
h_4 &=&  c_{32}^2c_{42} + c_{34}c_{42}^2 + c_{32}c_{42}c_{44} + c_{42}c_{44}^2
  + c_{22}c_{32} + c_{12}c_{42} + c_{24}c_{42} + c_{22}c_{44}
\end{eqnarray*}
On the other hand, in terms of the coordinates $d_{ij}$ we have
$$
M_D \;=\; \begin{pmatrix}
1 & d_{11} & 0 & d_{12}\\
0 & d_{21} & 1 & d_{22}\\
0 & d_{31} & 0 & d_{32}\\
0 & d_{41} & 0 & d_{42}
\end{pmatrix}\quad\hbox{\rm and}\quad  N_D \;=\;
\begin{pmatrix}
k_1 & 0 & \ell_1 & d_{13}\\
k_2 & 0 & \ell_2 & d_{23}\\
k_3 & 1 & \ell_3 & d_{33}\\
k_4 & 0 & \ell_4 & d_{43}
\end{pmatrix}
$$
where
\begin{eqnarray*}
k_1 &=&  d_{11}^2 + d_{12}d_{21} + d_{13}d_{31} + d_{14}d_{41}\\
k_2 &=&  d_{11}d_{21} + d_{21}d_{22} + d_{23}d_{31} + d_{24}d_{41}\\
k_3 &=&  d_{11}d_{31} + d_{21}d_{32} + d_{31}d_{33} + d_{34}d_{41}\\
k_4 &=&  d_{11}d_{41} + d_{21}d_{42} + d_{31}d_{42} + d_{41}d_{44}\\
\ell_1 &=&  d_{11}d_{12} + d_{12}d_{22} + d_{13}d_{32} + d_{14}d_{42}\\
\ell_2 &=&  d_{12}d_{21} + d_{22}^2 + d_{23}d_{32} + d_{24}d_{42}\\
\ell_3 &=&  d_{12}d_{31} + d_{22}d_{32} + d_{32}d_{33} + d_{34}d_{42}\\
\ell_4 &=&  d_{12}d_{41} + d_{22}d_{42} + d_{32}d_{43} + d_{42}d_{44}
\end{eqnarray*}
Now it is easy to compute the matrices $\widetilde{D}=M_C^{-1}\cdot N_C$
and $\widetilde{C}=M_D^{-1}\cdot N_D$. In order to compare~$\widetilde{C}$ to~$C$
and~$\widetilde{D}$ to~$D$, we have to transform from one coordinate system
to the other. For instance, we can substitute the entry $d_{ij}$ of~$D$ by
the $(i,j)$-entry of~$\widetilde{D}$ (which is a polynomial in the $c_{kl}$).
Upon performing this substitution in~$\widetilde{C}$, the result
should equal~$C$ modulo the ideal $I(\BO)$. Using \cocoa\ (cf.~\cite{CoCoA}),
it is straightforward to check that this is indeed the case.

On the side, we find that the open set $\mathbb{B}_{\mathcal{O},\mathcal{O}'}$
is defined as a subscheme of~$\mathbb{B}_{\mathcal{O}}$ by the non-vanishing of
$\det(M_C)=c_{24}c_{42}^2-c_{22}c_{42}c_{44}-c_{22}^2$, and the open set
$\mathbb{B}_{\mathcal{O}',\mathcal{O}}$ is defined as a subscheme 
of~$\mathbb{B}_{\mathcal{O}'}$ by the non-vanishing
of $\det(M_D)=d_{31}d_{42}-d_{32}d_{41}$.
\end{example}

The theorem allows us to answer to above question about the
relation between the matrices~$C_I$ and~$D_I$ representing~$I$
in the two $\mathcal{P}$-basis schemes. The explicit answer is given in the 
following corollary.

\begin{corollary}\label{basechangeforI}
Let~$I\subset P$ be an ideal such that both~$\mathcal{P}$ and~$\mathcal{P}'$
are bases modulo~$I$. Compute $M_{C_I},M_{D_I}\in \Mat_\mu(K)$,
$N_{C_I}\in\Mat_{\mu\alpha'}(K)$, and $N_{D_I}\in \Mat_{\mu,\alpha}(K)$ such that
$$
\mathcal{P}=\mathcal{P}'\cdot M_{D_I} \qquad
\mathcal{P}'=\mathcal{P}\cdot M_{C_I} \qquad
b\mathcal{P}=\mathcal{P}'\cdot N_{D_I} \qquad
b\mathcal{P}'=\mathcal{P}\cdot N_{C_I}
$$
hold in~$P/I$. Then $M_{C_I}$ and $M_{D_I}$ are invertible and we
have $C_I=M_{D_I}^{-1}\cdot N_{D_I}$ as well as $D_I=M_{C_I}^{-1}\cdot N_{C_I}$.
\end{corollary}

\begin{proof}
It suffices to substitute the coordinate tuples representing~$I$ in~$\BP$
and $\mathbb{B}_{\mathcal{P}'}$ in the matrix equalities given in part~(b) of
the theorem.
\end{proof}

A slight change in the point of view enables us to
determine the relation between the coefficients of the border bases
of two ideals which differ only by a linear change of coordinates.

\begin{proposition}\label{lc_represent}
Let $I\subset P$ be an ideal such that $\mathcal{P}$ is a basis modulo~$I$,
and let $\phi: P\longrightarrow P$ be a linear change of coordinates.
Write $\phi^{-1}(\mathcal{P})\equiv \mathcal{P}\cdot M_\phi\; (\mod I)$ with
a matrix $M_\phi\in\Mat_\mu(K)$, and write $\phi^{-1}(b\mathcal{P})\equiv\mathcal{P}\cdot
N_\phi\; (\mod I)$ with a matrix $N_\phi\in\Mat_{\mu,\alpha}(K)$. 
Then the following conditions are equivalent.
\begin{enumerate}
\item[(a)] The set~$\mathcal{P}$ is a basis modulo $\phi(I)$.
\item[(b)] The matrix~$M_\phi$ is invertible.
\end{enumerate}
In this case the ideal~$\phi(I)$ is represented 
by $C_{\phi(I)}=M_\phi^{-1}\cdot N_\phi$ in~$\BP$.
\end{proposition}

\begin{proof} 
By applying~$\phi$ to the congruence
$\phi^{-1}(\mathcal{P})\equiv \mathcal{P}\cdot M_\phi\; (\mod I)$,
we obtain $\mathcal{P}\equiv \phi(\mathcal{P})\cdot M_\phi\;
(\mod\, \phi(I))$. Using the fact that $\phi(\mathcal{P})$ 
is a basis modulo~$\phi(I)$, we see that~(a) and~(b) are
equivalent. Now we apply~$\phi$ to the second congruence 
in the proposition and get
$b\mathcal{P}\equiv \phi(\mathcal{P})\cdot N_\phi\;
(\mod\,\phi(I))$. By combining this with the previous
congruence, we obtain $b\mathcal{P} \equiv \mathcal{P}\cdot 
M^{-1}_\phi\cdot N_\phi\; (\mod\,\phi(I))$.
Thus~$\phi(I)$ is represented by $M_\phi^{-1}\cdot N_\phi$
in~$\BP$.
\end{proof}

At this point we can clarify the precise meaning of the idea
that a {\it generic}\/ linear change of coordinates should preserve
the property that~$I$ has an $\mathcal{O}$-basis and that
we should get a flat family in this way.
For this purpose, we introduce new indeterminates $a_{ij}$
for $i=1,\dots,n$ and $j=0,\dots,n$, and we let 
$A=(a_{ij})$. The $K$-algebra homomorphism
$$
\phi_A: K[a_{ij}][x_1,\dots,x_n] \longrightarrow K[a_{ij}][x_1,\dots,x_n]
$$
defined by $x_i \mapsto a_{i0}+a_{i1}x_1 + \cdots + a_{in}x_n$ is called
the {\bf generic linear change of coordinates}.
Letting $\widehat{\mathcal{A}}=(a_{ij})_{i,j=1,\dots,n}$, we see that
the set of linear changes of coordinates is the open subscheme~$\mathbb{L}$
of~$\mathbb{A}^{n(n+1)}$ defined by the non-vanishing 
of $\Delta=\det(\widehat{\mathcal{A}})$. We say that~$\mathbb{L}$ is the
{\bf scheme of linear coordinate changes}.
The coordinate ring of~$\mathbb{L}$ is $K[a_{ij}]_\Delta$.
Given a linear change of coordinates $\phi: P\longrightarrow P$
such that $\phi(x_i)=\alpha_{i0} +\alpha_{i1}x_1+\cdots +\alpha_{in}x_n$
with $\alpha_{ij}\in K$, we shall say that the matrix $\mathcal{A}=(\alpha_{ij})$,
or the point in~$\mathbb{L}$ whose coordinates are the entries 
of~$\mathcal{A}$, represent~$\phi$ in~$\mathbb{L}$.

\begin{proposition}\label{linchange}
Let $I\subset P$ be an ideal such that $\mathcal{P}$ is a basis modulo~$I$,
let~$\phi_A$ be the generic linear change of coordinates, let
$\bar I = I\cdot K[a_{ij}]_\Delta[x_1,\dots,x_n]$, and
let~$W$ be the Zariski-open subset of~$\mathbb{L}$ 
consisting of all points representing linear changes of 
coordinates $\phi: P\longrightarrow P$ 
such that the matrix~$M_\phi$ is invertible.

\begin{enumerate}
\item[(a)] The open set $W$ is non-empty.

\item[(b)] Let $M_{\phi_A} \in\Mat_\mu(K[a_{ij}]_\Delta)$ be such that
$\phi^{-1}_A(\mathcal{P})\equiv\mathcal{P} \cdot M_{\phi_A} (\mod \bar I)$,
and let $\Lambda=\det(M_{\phi_A})$. Then
the coordinate ring of~$W$ is~$K[a_{ij}]_{\Delta\cdot\Lambda}$.

\item[(c)] The entries of the matrix $M_{\phi_A}^{-1}\cdot N_{\phi_A}$ 
are contained in~$K[a_{ij}]_{\Delta\cdot\Lambda}$. 
Hence there exists a well-defined $K$-algebra 
homomorphism $\psi: B_{\mathcal{P}} \longrightarrow 
K[a_{ij}]_{\Delta\cdot\Lambda}$ which satisfies 
$\psi(c_{ij})= (M_{\phi_A}^{-1}\cdot N_{\phi_A})_{i,j}$ for all $i,j$.

\item[(d)] The map~$\psi$ induces a morphism  $\Phi: W \longrightarrow \BP$
of affine schemes which is defined as follows. 
If~$p\in W$ (resp.\ the corresponding matrix~$\mathcal{A}$)
represents a linear change of coordinates~$\phi: P\longrightarrow P$,
then~$\Phi(p)$ is represented by $C_{\phi(I)}=M_\phi^{-1}
\cdot N_{\phi}$. 

\item[(e)] The map~$\psi$ induces a flat family
$\Psi: K[a_{ij}]_{\Delta\cdot\Lambda}\longrightarrow 
U_{\mathcal{P}} \otimes_{B_{\mathcal{P}}} K[a_{ij}]_{\Delta\cdot\Lambda}$.
\end{enumerate}
\end{proposition}

\begin{proof} 
For the proof of~(a), it suffices to observe that~$W$ contains 
the point corresponding to the identity map. To prove~(b),
we apply the preceding proposition. It says that, for a point
$(\alpha_{ij})\in \mathbb{L}$ representing a linear change of
coordinates~$\phi$, the set~$\mathcal{P}$ is a basis modulo~$\phi(I)$ 
if and only if $\det(M_\phi)=\Lambda(\alpha_{ij})\ne 0$.

Next we show~(c). The fact that the entries of $M_{\phi_A}^{-1}$ are 
in~$K[a_{ij}]_{\Delta\cdot\Lambda}$ follows from the definition.
From $\phi_A^{-1}(b\mathcal{P})\equiv \mathcal{P}\cdot N_{\phi_A}\;
(\mod \bar I)$ we see that the entries of~$N_{\phi_A}$ 
are polynomials in the rational functions $\frac{a_{ij}}{\Delta}$.
These observations show that~$\psi(c_{ij})\in K[a_{ij}]_{\Delta\cdot\Lambda}$ 
for all~$i,j$.

Claim~(d) follows immediately from~(c), and~(e)
is a consequence of the flatness of the universal family
(see Proposition~\ref{PseudoUF}) by applying a base change 
with~$K[a_{ij}]_{\Delta\cdot\Lambda}$.
\end{proof}

The following example illustrates the proposition.

\begin{example}
Let $I$ be the ideal in $K[x_1,x_2]$
generated by $\{x_1^2-1,\, x_2^2-1\}$. 
Then $\mathcal{O} = (1, x_1, x_2, x_1x_2)$ is a basis modulo~$I$. 
The generic linear change of coordinates~$\phi_A$ is given by
\begin{eqnarray*}
\phi_A(x_1) &=& a_{10} + a_{11}\, x_1 + a_{12}\,x_2\\
\phi_A(x_2) &=& a_{20} + a_{21}\, x_1 + a_{22}\, x_2
\end{eqnarray*}
We have $\Delta = a_{11}a_{22}-a_{12}a_{21}$. The inverse of~$\phi_A$ satisfies
$$
\begin{matrix}
\phi_A^{-1}(x_1) &=& \tfrac{1}{\Delta}\, (a_{22}\, x_1 - a_{12}\, x_2) - b_{10}\\
\phi_A^{-1}(x_2) &=& \tfrac{1}{\Delta}\, (-a_{21}\, x_1 + a_{11}\, x_2) - b_{20}\\
\end{matrix}\qquad\hbox{\rm with}\quad
\binom{b_{10}}{b_{20}}=\widehat{A}^{-1}\cdot \binom{a_{10}}{a_{20}}
$$
Using the fact that~$\mathcal{O}$ is a basis modulo~$I$,
we write $\phi_A^{-1}(\mathcal{O}) = \mathcal{O}\cdot M_{\phi_A}
\; (\mod \bar I)$ where
$$
M_{\phi_A} =
\begin{pmatrix}
1 & -b_{10}        & -b_{20}        & b_{10}b_{20}-(a_{11}a_{12}+a_{21}a_{22})/\Delta^2 \cr
0 & a_{22}/\Delta  & -a_{21}/\Delta & (a_{21}b_{10}-a_{22}b_{20})/\Delta  \cr
0 & -a_{12}/\Delta & a_{11}/\Delta  & (a_{12}b_{20}-a_{11}b_{10})/\Delta \cr
0 & 0      & 0      & (a_{11}a_{22}+a_{12}a_{21})/\Delta^2
\end{pmatrix}
$$
Thus we have $\Lambda = \det(M_{\phi_A})= (a_{11}a_{22} +a_{12}a_{21})/\Delta^3$. 

Now we consider a specific $K$-algebra homomorphism
$\phi:P\longrightarrow P$ given by $x_i\mapsto
\alpha_{i0}+\alpha_{i1}x_1+\alpha_{i2}x_2$
with $\alpha_{ij}\in K$. The condition that~$\phi$
is a linear change of coordinates is expressed 
by~$\Delta(\alpha_{ij})\ne 0$. The additional condition that 
$M_\phi$ is invertible is then expressed by $\Lambda(\alpha_{ij})\ne 0$,
because we have $M_\phi= M_{\phi_A} \vert_{a_{ij}\mapsto 
\alpha_{ij}}$.

For instance, let $\phi: K[x_1,x_2]\longrightarrow K[x_1,x_2]$
be given by $\phi(x_1)=x_1+x_2$ and $\phi(x_2)=x_1-x_2$, i.e.\ let
$\alpha_{10} = \alpha_{20} = 0$, 
$\alpha_{11} = \alpha_{22} = \alpha_{21} = 1$, and
$\alpha_{12} = -1$. Then $\Delta(\alpha_{ij})\ne 0$ shows that~$\phi$ is invertible.
Now $\Lambda(\alpha_{ij}) = 0$ implies that~$M_\phi$ 
is not invertible. Hence~$\mathcal{O}$ is not a basis modulo~$\phi(I)$.
In fact, if we perform the linear change, we see that~$\phi(I)$ is 
generated by $\{(x_1-x_2)^2-1,\, (x_1+x_2)^2-1\}$, and therefore
by $\{x_1^2+x_2^2-1,\, x_1x_2\}$. Since $x_1x_2\in\phi(I)$, it is clear
that~$\mathcal{O}$ is not a basis modulo $\phi(I)$.
\end{example}

The existence of a flat family of ideals defined by
linear changes of coordinates distinguishes border bases
from Gr\"obner bases in the following sense.

\begin{corollary}\label{no-gin-cor}
Let $\mathcal{O}$ be an order ideal in~$\mathbb{T}^n$, and let 
$I\subset P$ be an ideal which has an $\mathcal{O}$-border basis.
Then, for a generically chosen linear change of coordinates~$\phi$,
the ideal~$\phi(I)$ has again an $\mathcal{O}$-border basis.
\end{corollary}

\begin{proof}
This follows from Proposition~\ref{lc_represent} and
the fact that $W\ne \emptyset$ by the preceding proposition.
\end{proof}

Notice that the property described in this corollary 
differs markedly from Gr\"obner basis theory
where a generically chosen linear change of coordinates
entails in general a new leading term ideal, and therefore also
a new order ideal $\mathcal{O}_\sigma(I)=\mathbb{T}^n \setminus
\LT_\sigma(I)$. In border basis theory there is {\it no gin!}

Finally, we point out two particular situations in which
the claims of the preceding two propositions simplify
substantially.

\begin{example}
Let us consider the set of all translations, i.e.\ of the
linear changes of coordinates with $\widehat{A}=I_n$.
They are all invertible  and their inverses are also translations. 
If ~$\phi$ is a translation and we order the elements of~$\mathcal{O}$ 
in increasing degree, then~$M_\phi$ is an upper triangular matrix 
having all entries on the main diagonal equal to~$1$. Therefore the 
matrix~$M_\phi$ is invertible for every translation~$\phi$. 
Hence, given an ideal $I\subset P$ such that $\mathcal{O}$ is a basis modulo~$I$,
the order ideal~$\mathcal{O}$ is also a basis modulo $\phi(I)$.
\end{example}

\begin{example}
Consider the order ideal $\mathcal{O}=\{1,x_1,\dots,x_n\}$, and
let $\phi:P\longrightarrow P$ be a linear change of coordinates.
We write $\phi(x_i)=a_{i0}+a_{i1}x_1 +\cdots + a_{in}x_n$
with $a_{ij}\in K$.
Given an ideal~$I$ which has an $\mathcal{O}$-border basis, 
the matrix~$M_{\phi_A}$ is defined by $\phi^{-1}(\mathcal{O})
\equiv\mathcal{O}\cdot M_{\phi_A}\; (\mod I)$. Here we get 
$$
M_{\phi_A}= \begin{pmatrix}
1 & -b_{10}\;\cdots\;-b_{n0}\\
0 &  \\
\vdots & (\widehat{A}^{-1})^{\rm tr} \\
0 & \end{pmatrix} 
\quad\hbox{\rm where}\quad 
\begin{pmatrix}b_{10} \\ \vdots \\ b_{n0}\end{pmatrix} =
\widehat{A}^{-1} \cdot 
\begin{pmatrix}a_{10} \\ \vdots \\ a_{n0}\end{pmatrix}
$$
Hence we have $W=\mathbb{L}$ in Proposition~\ref{linchange}.
In other words, if $\mathcal{O}$ is a basis modulo~$I$ 
and $\phi: P \longrightarrow P$ is a linear change of coordinates,
then $\mathcal{O}$ is also a basis modulo~$\phi(I)$.
\end{example}

\bigskip
%%%%%%%%%%%%%%%%%%%%%%%%%%%%%%%%%%%%%%%%%%%%%%%%%%%%%%%
%
% Section 3: The Principal Component of the BB Scheme
%
%%%%%%%%%%%%%%%%%%%%%%%%%%%%%%%%%%%%%%%%%%%%%%%%%%%%%%%

\section{The Principal Component of the Border Basis Scheme}
\label{The Principal Component of the Border Basis Scheme}
\bigskip

As mentioned in the introduction, our next goal is to study
the principal component of the border basis scheme.
Since we do not need the very general setting of the previous
section, we shall concentrate on the classical border basis scheme.

In the following we continue to work over the
polynomial ring $P=K[x_1,\dots,x_n]$ over a field~$K$,
we let~$\mathcal{O}=\{t_1,\dots,t_\mu\}$ be an order ideal
of terms in~$\mathbb{T}^n$, and we let
$\partial\mathcal{O}=\{b_1,\dots,b_\nu\}$ be its border.
Moreover, we denote the algebraic closure of~$K$ by~$\overline{K}$,
and we let $\overline{P}=\overline{K}[x_1,\dots,x_n]$.

\begin{definition}\label{defprincomp}
For each 0-dimensional ideal $I\subset \overline{P}$ having an 
$\mathcal{O}$-border basis, let $\beta(I)$ be the corresponding point
of~$\BO\times_{\Spec(K)}\Spec(\overline{K})$. Then the 
closed subscheme~$\CO$ of~$\BO$ such that
$\CO\times_{\Spec(K)}\Spec(\overline{K})$
is the closure of the set of all points
$\beta(I_{\mathbb{X}})$, where $\mathbb{X}\subseteq \mathbb{A}^n(\overline{K})$
is a reduced scheme of length~$\mu$ and $I_{\mathbb{X}}\subset \overline{P}$ is 
its vanishing ideal, is called the {\bf principal component} of~$\BO$.
\end{definition}

It is known that the {\it radical component}\/ of the Hilbert scheme
is irreducible (see~\cite{MS}, 18.32).
Since~$\CO$ is a Zariski-open subset of the radical component,
it follows that~$\CO$ is an irreducible component of~$\BO$,
so that its name is justified.
This result is also shown in Theorem~\ref{princicomp} below.

As promised in the introduction, we
will construct explicit equations defining~$\CO$.
Our method is inspired by suggestions in~\cite{Ha1}, p.~213
and~\cite{Ha2}, Sect.~2.1. We use additional indeterminates~$y_j^{(i)}$
for $i=1,\dots,\mu$ and $j=1,\dots,n$ and we group them into
tuples $\yi=(y_1^{(i)},\dots,y_n^{(i)})$. The indeterminates
in~$\yi$ should be thought of as representing the coordinates
of the $i^{\rm th}$ point of~$\mathbb{X}$.

\begin{definition}
Let $Q=K[\y1,\dots,\ymu]$. We define the following
polynomials in~$Q$.

\begin{enumerate}
\item[(a)] Let $\Delta_{\mathcal{O}}= \det (t_j(\yi))_{i,j=1,\dots,\mu}$
where $t_j(\yi)$ denotes the result of the substitutions
$x_k\mapsto y_k^{(i)}$ in~$t_j$.

\item[(b)] For $i=1,\dots,\mu$ and $j=1,\dots,\nu$,
we let
$$
\Delta_{ij}=\det (t_1(\yd) \mid \cdots \mid
b_j(\yd) \mid \cdots \mid t_\mu(\yd))
$$
Here $t_k(\yd)$ denotes the $k^{\,\rm th}$ column of the
matrix $\Delta_{\mathcal{O}}$, i.e.\ the column
$(t_k(\y1),\dots,t_k(\ymu))^{\rm tr}$. Thus~$\Delta_{ij}$
is the determinant of the matrix where the $i^{\,\rm th}$ column
of $\Delta_{\mathcal{O}}$ has been replaced by $b_j(\yd)$.)
\end{enumerate}
\end{definition}

This definition can be motivated as follows.

\begin{remark}\label{pij_cij_corresp}
Notice that $\Delta_{\mathcal{O}}\ne 0$, since
each row contains different indeterminates.

\begin{enumerate}
\item[(a)]
In the quotient field of~$Q$,
consider the system of linear equations
\begin{eqnarray*}
t_1(\y1)\,z_1 + \cdots + t_\mu(\y1)\, z_\mu & = & b_j(\y1) \\
\vdots\qquad\qquad\qquad &&\quad\;\vdots \\
t_1(\ymu)\,z_1 + \cdots + t_\mu(\ymu)\, z_\mu & = & b_j(\ymu) \\
\end{eqnarray*}
By Cramer's Rule, its solution is given
by $\textstyle{\frac{1}{\Delta_{\mathcal{O}}}}
\cdot (\Delta_{1j},\dots,\Delta_{\mu j})$.

\item[(b)] Given a set of points $\mathbb{X}=
\{\mathbf{p}_1,\dots,\mathbf{p}_\mu\}$ whose vanishing ideal~$I_{\mathbb{X}}$
has an $\mathcal{O}$-border basis, we can substitute the
coordinates $p_{ij}$ of the points
$\mathbf{p}_i=(p_{i1},\dots,p_{in})$ for the indeterminates $y_j^{(i)}$
in the systems of linear equations above.
The solutions $(\gamma_{ij})$ of the resulting systems
are precisely the coefficients of the border basis
$G=\{g_1,\dots,g_\nu\}$ of the ideal~$I_{\mathbb{X}}$. Here we
have $g_j=b_j -\sum_i \gamma_{ij}t_i$.

\end{enumerate}
\end{remark}

The main result of this subsection is that the following
ring is isomorphic to the affine coordinate ring of the
principal component of~$\BO$.

\begin{notation}
Let $K(\y1,\dots,\ymu)$ be the quotient field of~$Q$.

\begin{enumerate}
\item[(a)] Let $C_{\mathcal{O}}$ be the $K$-subalgebra
of $K(\y1,\dots,\ymu)$ generated by the elements
$\Delta_{ij}/\Delta_{\mathcal{O}}$ with $i\in\{1,\dots,\mu\}$
and $j\in\{1,\dots,\nu\}$.

\item[(b)]  For $i\in\{1,\dots,\mu\}$ and $j\in\{1,\dots,\nu\}$,
let $c_{ij}$ be new indeterminates.
We define a surjective $K$-algebra homomorphism
$\Phi: K[c_{ij}] \longrightarrow C_{\mathcal{O}}$
by letting $\Phi(c_{ij})=\Delta_{ij}/\Delta_{\mathcal{O}}$.
\end{enumerate}
\end{notation}

\begin{lemma}\label{IBinker}
The defining ideal of~$\BO$ is contained in the kernel
of~$\Phi$. Consequently, the map~$\Phi$
induces a surjective $K$-algebra homomorphism
$B_{\mathcal{O}} \longrightarrow C_{\mathcal{O}}$
and a closed immersion $\Spec(C_{\mathcal{O}}) \hookrightarrow
\BO$.
\end{lemma}

\begin{proof}
The ideal~$I(\BO)$ defining~$\BO$ is generated by
the entries of the commutators $\mathcal{A}_k\mathcal{A}_\ell
-\mathcal{A}_\ell \mathcal{A}_k$ of the formal multiplication
matrices of the generic $\mathcal{O}$-border basis.
Thus we have to show that the matrices $\Phi(\mathcal{A}_k)$
commute, where $\Phi(\mathcal{A}_k)$ is obtained
by applying~$\Phi$ to the entries of the matrix~$\mathcal{A}_k$.

The $j^{\,\rm th}$ column of $\Phi(\mathcal{A}_k)$ is the
solution of the system of linear equations
$$
(t_1(\yd) \mid \cdots \mid t_\mu(\yd)) \cdot (z_1,\dots,z_\mu)^{\rm tr}
= (x_kt_j)(\yd)
$$
Therefore we get the following equalities $(\ast)$:
\begin{eqnarray*}
(t_1(\yd) \mid \cdots \mid t_\mu(\yd))\cdot \Phi
(\mathcal{A}_k)
& = & ((x_kt_1)(\yd) \mid \cdots \mid (x_k t_\mu)(\yd)) \\
&=& {\rm diag}(y_k^{(1)},\dots,y_k^{(\mu)}) \cdot (t_1(\yd) \mid
\cdots \mid t_\mu(\yd))
\end{eqnarray*}
Since diagonal matrices commute, it follows that
$$
(t_1(\yd) \mid \cdots \mid t_\mu(\yd))\cdot \Phi(\mathcal{A}_k)
\cdot \Phi(\mathcal{A}_\ell) =
(t_1(\yd) \mid \cdots \mid t_\mu(\yd))\cdot \Phi(\mathcal{A}_\ell)
\cdot \Phi(\mathcal{A}_k)
$$
and the fact that the matrix $(t_1(\yd) \mid \cdots \mid t_\mu(\yd))$
is invertible over the quotient field of~$Q$ implies the claim.
\end{proof}

In fact, the image of the closed immersion we just found is exactly
the principal component of the border basis scheme, as our next 
theorem shows.

\begin{theorem}{\rm (The Coordinate Ring of the Principal Component)}
\label{princicomp}\\
Let $\Phi: K[c_{ij}] \longrightarrow C_{\mathcal{O}}$
be the surjective $K$-algebra homomorphism defined by 
$\Phi(c_{ij}) = \Delta_{ij}/\Delta_{\mathcal{O}}$.

\begin{enumerate}
\item[(a)] The ideal~$\ker(\Phi)$ is the vanishing ideal of the
principal component~$\CO$ of the border basis scheme~$\BO$.
In particular, the principal component~$\CO$
of~$\BO$ is the closure of the image of the morphism
$\Spec(C_{\mathcal{O}}) \hookrightarrow\BO$.

\item[(b)] The scheme~$\CO$ is irreducible.
\end{enumerate}
\end{theorem}

\begin{proof}
Let $\overline{K}$ be the algebraic closure of~$K$.
Since the base change $K \subseteq \overline{K}$ is faithfully
flat, it suffices to prove that the ideal
$\ker(\Phi)\cdot \overline{K}[c_{ij}]$
is the vanishing ideal of the scheme $\CO \times_{\Spec(K)}
\Spec(\overline{K})$.
In other words, we may (and shall)
assume that~$K$ is algebraically closed.

First we show that the map $\Spec(\Phi)$ yields a
bijection between the closed points
of~$\Spec(C_{\mathcal{O}})$ and the closed points of
$\CO$. A closed point $(p_{ij})$ of $\Spec(C_{\mathcal{O}})$
corresponds to a maximal ideal $\mathfrak{m}=\langle
y_j^{(i)}- p_{ij} \rangle_{i=1,\dots,\mu; j=1,\dots,n}$
of~$K[\y1,\dots,\ymu]$ which does not contain~$\Delta_{\mathcal{O}}$. 
Thus~$\mathfrak{m}$ defines
also a maximal ideal in the localization
$K[\y1,\dots,\ymu]_{\Delta_{\mathcal{O}}}$ and by intersecting
it with~$C_{\mathcal{O}}$ we obtain a maximal ideal
$\mathfrak{m}_{\mathcal{O}}$ of~$C_{\mathcal{O}}$.

Let us examine this maximal ideal~$\mathfrak{m}_{\mathcal{O}}$
more closely. The elements $p_{ij}\in K$ define
a set of points $\mathbb{X}=\{\mathbf{p}_1,\dots,\mathbf{p}_\mu\}$ where
$\mathbf{p}_i=(p_{i1},\dots,p_{in})$. The image of~$\Delta_\mathcal{O}$
in $K[\y1,\dots,\ymu]/\mathfrak{m}$ is the determinant
of~$(t_j(\mathbf{p}_i))_{i,j}$. Since we assume that this determinant
is non-zero, the set~$\mathbb{X}$ consists of pairwise distinct points.
Moreover, it follows that the ideal~$I_{\mathbb{X}}$ has an
$\mathcal{O}$-border basis.

The systems of linear equations
$(t_1(\mathbf{p}_i) \mid \cdots \mid t_\mu(\mathbf{p}_i))
\cdot (z_1,\dots,z_\mu)^{\rm tr} = (b_j(\mathbf{p}_i))$
have unique solutions $(\gamma_{1j},\dots,\gamma_{\mu j})\in K^\mu$.
Then the  corresponding $\mathcal{O}$-border prebasis
$G=\{g_1,\dots,g_\nu\}$ with $g_j=b_j -\sum_{i=1}^\mu
\gamma_{ij} t_i$ is the $\mathcal{O}$-border basis
of~$I_{\mathbb{X}}$ and the maximal ideal
$\langle c_{ij}-\gamma_{ij}\rangle_{i,j}$ of~$K[c_{ij}]$
defines a closed point of~$\CO$.
By Remark~\ref{pij_cij_corresp}.b, the maximal ideal
$\langle c_{ij}-\gamma_{ij}\rangle_{i,j}$ is precisely
the preimage of~$\mathfrak{m}$ under~$\Phi$.

Conversely, let us start with a closed point $(\gamma_{ij})$ 
of~$\BO$ corresponding to the $\mathcal{O}$-border basis of a 
radical ideal~$I$. Since the base field is algebraically closed,
the ideal~$I$ is the vanishing ideal of a set of~$\mu$
points $\mathbb{X}=\{(p_{i1},\dots,p_{in}) \mid i=1,\dots,\mu\}$ 
in~$\mathbb{A}^n$. From what we just showed
it follows that the maximal ideal $\langle c_{ij}-\gamma_{ij}\rangle_{i,j}$ 
is the preimage of the maximal ideal  
$\langle y_j^{(i)}- p_{ij} \rangle_{i=1,\dots,\mu; j=1,\dots,n}$
under~$\Phi$.

Next we let $R=K[\Delta_{ij} \mid i\in\{1,\dots,\mu\},
j\in\{1,\dots,\nu\}]_{\Delta_{\mathcal{O}}}$ and
consider the canonical injective $K$-algebra homomorphism
$C_{\mathcal{O}} \hookrightarrow R$. The intersection of the
maximal ideals of~$R$ is~$(0)$. The preimages of these
maximal ideals in~$C_{\mathcal{O}}$ are exactly the maximal
ideals~$\mathfrak{m}_{\mathcal{O}}$ constructed above.
Hence the intersection of the maximal ideals
$\mathfrak{m}_{\mathcal{O}}$ is $(0)$. From what we have shown,
it follows that the intersection of the maximal ideals
$\langle c_{ij}-\gamma_{ij}\rangle_{i,j}$ of~$K[c_{ij}]$
is the kernel of~$\Phi$.

On the other hand, the general form of Hilbert's Nullstellensatz 
(see~\cite{KR1}, 2.6.17) implies that the intersection of the
maximal ideals $\langle c_{ij}-\gamma_{ij}\rangle_{i,j}$ of~$K[c_{ij}]$
is the vanishing ideal of the closure of the corresponding set of
points, i.e.\ the vanishing ideal of the principal component~$\CO$
of~$\BO$.

To prove~(b), we note that the ring $C_{\mathcal{O}}$ is a subalgebra 
of $K(\y1,\dots,\ymu)$. Hence it is an integral domain. Therefore the 
ideal $\ker(\Phi)$ is a prime ideal. 
\end{proof}

In view of the preceding theorem it will prove useful
to study the $K$-algebra $C_{\mathcal{O}}$ in more detail.
Our next proposition shows that it contains the
following elements.

\begin{notation}
Let $L=\{s_1,\dots,s_\mu\}$ be a set of~$\mu$ distinct
terms contained in $\mathbb{T}^n$.
Then we set $\Delta_L=\det (s_1(\yi) \mid \cdots \mid s_\mu(\yi))
\in Q$.
\end{notation}

\begin{proposition}\label{COgens}
For every $L=\{s_1,\dots,s_\mu\} \subset \mathbb{T}^n$,
we have $\Delta_L/\Delta_{\mathcal{O}} \in C_{\mathcal{O}}$.
\end{proposition}

\begin{proof}
If $L=\mathcal{O}$, we have $\Delta_L/\Delta_{\mathcal{O}}=1
\in C_{\mathcal{O}}$. Next we show that
$\Delta_{L_j}/\Delta_{\mathcal{O}} \in C_{\mathcal{O}}$
for $L_j=(t_1,\dots,t_{j-1},s,t_{j+1},\dots,t_\mu)$ with
$j\in\{1,\dots,\mu\}$ and a term
$s\in \mathbb{T}^n\setminus \mathcal{O}$. We write $s=t'b_j$ with
$t'\in \mathbb{T}^n$, $j\in\{1,\dots,\nu\}$ and
we prove the claim by induction on $\deg(t')$.

If $\deg(t')=0$, the term~$s$ is a border term and~$\Delta_L$
is one of the elements $\Delta_{ij}$. Now let $\deg(t')>0$ and write
$t'=x_k\,t''$ with $k\in\{1,\dots,n\}$ and $t''\in\mathbb{T}^n$.
By Cramer's rule, we have
\begin{gather*}
(t_1(\yd) \mid \cdots \mid t_\mu(\yd)) \cdot (\Delta_{L_1}/
\Delta_{\mathcal{O}},\dots,
\Delta_{L_\mu}/\Delta_{\mathcal{O}})^{\rm tr} \; =\; ((x_kt''b_j)(\yd)
\qquad\qquad \\ \qquad\qquad\qquad\qquad\qquad\qquad\qquad\qquad
=\; {\rm diag}(y_k^{(1)},\dots,y_k^{(\mu)}) \cdot ((t''b_j)(\yd))
\end{gather*}
The inductive hypothesis implies that there are elements
$\widetilde{\Delta}_1/\Delta_{\mathcal{O}}, \dots,
\widetilde{\Delta}_\mu/\Delta_{\mathcal{O}} \in C_{\mathcal{O}}$
such that
$$
(t_1(\yd) \mid \cdots \mid t_\mu(\yd)) \cdot
(\widetilde{\Delta}_1/\Delta_{\mathcal{O}}, \dots,
\widetilde{\Delta}_\mu/\Delta_{\mathcal{O}})^{\rm tr} =
((t''b_j)(\yd))
$$
Using the equality $(\ast)$ from the proof of Lemma~\ref{IBinker},
we get
\begin{gather*}
(t_1(\yd) \mid \cdots \mid t_\mu(\yd)) \cdot (\Delta_{L_1}/
\Delta_{\mathcal{O}},
\dots,\Delta_{L_\mu}/\Delta_{\mathcal{O}})^{\rm tr} \;=
\qquad\qquad\qquad\qquad\qquad\qquad \\
\qquad\qquad =\;
{\rm diag}(y_k^{(1)},\dots,y_k^{(\mu)}) \cdot
(t_1(\yd) \mid \cdots \mid t_\mu(\yd)) \cdot
(\widetilde{\Delta}_1/\Delta_{\mathcal{O}}, \dots,
\widetilde{\Delta}_\mu/\Delta_{\mathcal{O}})^{\rm tr} \\
\qquad\qquad =\;
(t_1(\yd) \mid \cdots \mid t_\mu(\yd)) \cdot \widetilde{\Phi}(\mathcal{A}_k) \cdot
(\widetilde{\Delta}_1/\Delta_{\mathcal{O}}, \dots,
\widetilde{\Delta}_\mu/\Delta_{\mathcal{O}})^{\rm tr}
\qquad\qquad\quad
\end{gather*}
At this point we note that $(t_1(\yd) \mid \cdots \mid t_\mu(\yd))$
is an invertible matrix over the field $K(\y1,\dots,\ymu)$.
It follows that the tuple $(\Delta_{L_1}/\Delta_{\mathcal{O}},
\dots,\Delta_{L_\mu}/\Delta_{\mathcal{O}})$ is contained in
$(C_{\mathcal{O}})^\mu$.

Finally we turn to the general case. Let $L=\{t_{i_1},\dots,t_{i_k},
s_1,\dots,s_{\mu-k}\}$ with $i_1,\dots,i_k\in \{1,\dots,\mu\}$
and $s_j\in\mathbb{T}^n\setminus \mathcal{O}$.
Clearly, we may assume that the indices $i_1,\dots,i_k$ are pairwise
distinct. We proceed by downward induction on~$k$. The case $k=\mu-1$
has been treated above. For the induction step, let $\{i_{k+1},\dots,
i_{\mu-k}\}= \{1,\dots,\mu\}\setminus\{i_1,\dots,t_k\}$.
Now the claim follows from the Pl\"ucker relation
\begin{eqnarray*}
\Delta_L\cdot \Delta_{\mathcal{O}} &=& \det(t_{i_1}(\yd) \mid
\cdots \mid t_{i_k}(\yd) \mid s_1(\yd) \mid\cdots\mid
s_{\mu-k}(\yd))\\
&& \;\cdot\; \det(t_{i_1}(\yd) \mid\cdots\mid
t_{i_\mu}(\yd)) \\
&=& {\!\!\!\textstyle\sum\limits_{j=1}^{\mu-k}}
\pm \det(t_{i_1}(\yd) \mid \cdots\mid t_{i_{k+1}}(\yd) \mid
s_1(\yd) \mid \cdots\mid \widehat{s_j(\yd)} \mid \cdots\mid
s_{\mu-k}(\yd)) \\
&& \cdot\; \det(t_{i_1}\mid\cdots\mid
\widehat{t_{i_k}(\yd)} \mid \cdots\mid t_{i_{\mu-k}}(\yd)\mid
s_j(\yd))
\end{eqnarray*}
and the inductive hypothesis.
\end{proof}

In other words, this proposition says that
$C_{\mathcal{O}}=K[\Delta_L/\Delta_{\mathcal{O}} 
\mid L\subset\mathbb{T}^n,\, \#L=\mu]$. 
Therefore the ring~$C_{\mathcal{O}}$ agrees with
the one mentioned in~\cite{Ha1} and~\cite{Ha2}.
Restricting the number of algebra generators has an obvious
advantage: we can now write down an algorithm for computing
the defining equations of the principal component.
This makes it possible to check effectively whether a given
border basis scheme is irreducible.

\begin{proposition}\label{PCcomp}
Let $\mathcal{O}=\{t_1,\dots,t_\mu\}$ be an order ideal
in~$\mathbb{T}^n$ and let $c_{ij}$ be further indeterminates,
where $i=1,\dots,\mu$ and $j=1,\dots,\nu$.
The following instructions define an algorithm
which computes a system of generators of the defining ideal
in~$K[c_{ij}]$ of the principal component~$\CO$ of the
border basis scheme.

\begin{enumerate}
\item[(1)] Form the polynomial ring $Q=K[\y1,\dots\ymu]$
and compute the elements $\Delta_{\mathcal{O}}=
\det(t_j(\yi)$ and $\Delta_{ij}=\det (t_1(\yd) \mid \cdots \mid
b_j(\yd) \mid \cdots \mid t_\mu(\yd))$ for $i=1,\dots,\mu$
and $j=1,\dots,\nu$.

\item[(2)]  Form the polynomial ring $K[c_{ij},z]$ where~$z$
is a new indeterminate. Let~$I$ be the ideal generated by
$\Delta_{\mathcal{O}}\, z-1$ and the set of all
$\Delta_{\mathcal{O}}\, c_{ij} - \Delta_{ij}$
such that $i=1,\dots,\mu$ and $j=1,\dots,\nu$.

\item[(3)] Compute a set of generators $\{F_1,\dots, F_r\}$ 
of the elimination ideal $I\cap K[c_{ij}]$ and return it.

\end{enumerate}
\end{proposition}

\begin{proof}
This is a special case of a classical algorithm which computes explicit
representations of finitely generated subalgebras of
function fields (see for instance~\cite{KR2}, Tutorial 41).
\end{proof}

\begin{proposition}
Let $W\in\Mat_{m,n}(\mathbb{Z})$, and let~$P$ be graded by~$W$.
We equip $K[c_{ij}]$ with a $\mathbb{Z}^m$-grading given by a 
matrix~$\overline{W}$ such that $\deg_{\overline{W}}(c_{ij}) = 
\deg_W(b_j) - \deg_W(t_i)$.
Then the elimination ideal $I\cap K[c_{ij}]$ of the preceding proposition is  
homogeneous with respect to the grading given by~$\overline{W}$.
\end{proposition}

\begin{proof}
First we introduce a $\mathbb{Z}^m$-grading on $Q=K[\y1,\dots\ymu]$
given by a matrix $\widetilde{W}$
by letting $\deg_{\widetilde{W}}(y_j^{(i)})=\deg_W(x_i)$.
Thus the elements of the $j^{\rm th}$ column of the matrix~$t_j(\yi)$
are homogeneous of degree $\deg_W(t_j)$. Hence~$\Delta_{\mathcal{O}}$
is homogeneous of degree $\deg_W(t_1\cdots t_\mu)$.
Similarly, we see that $\Delta_{ij}$ is homogeneous of degree
$\deg_W(t_1\cdots t_\mu)-\deg_W(t_i)+\deg_W(b_j)$.
This shows that if we define $\deg_{\overline{W}}(c_{ij}) = \deg_W(b_j) - \deg_W(t_i)$
and $\deg_{\overline{W}}(z) = -\deg_W(t_1\cdots t_\mu)$, then~$I$ is a homogeneous ideal
in $K[c_{ij},z]$. Consequently, also $I\cap K[c_{ij}]$ is a homogeneous ideal  
in~$K[c_{ij}]$ with respect to the grading given by~$\overline{W}$.
\end{proof}

This result is in accordance with the fact that
the ideal~$I(\BO)$ is homogeneous with respect to
the same grading.

\begin{remark}
Suppose that the order ideal~$\mathcal{O}$ is a $\sigma$-cornercut
with respect to some term ordering~$\sigma$. This implies that
there exists a system of positive weights for $x_1,\dots,x_n$  
such that $\deg_W(b_j)> \deg_W(t_i)$ for $i=1,\dots,\mu$ and
$j=1,\dots,\nu$. By the proposition, the ideal
$I\cap K[c_{ij}]$ is homogeneous with respect to a  
positive grading on~$K[c_{ij}]$. This observation agrees with 
the fact that, in this case, the border basis scheme 
$\mathbb{B}_\mathcal{O}$ and the Gr\"obner basis scheme
$\mathbb{G}_{\mathcal{O},\sigma}$ are isomorphic (see~\cite{R}, 
Proposition 3.11), and the latter can be seen as a
weighted projective scheme (see~\cite{R}, Theorem 2.8).
\end{remark}

\bigskip
%%%%%%%%%%%%%%%%%%%%%%%%%%%%%%%%%%%%%%%%%%%%%%%%%%%%%%%
%
% Section 4: Local Parameters at the Radical Points
%
%%%%%%%%%%%%%%%%%%%%%%%%%%%%%%%%%%%%%%%%%%%%%%%%%%%%%%%

\section{Local Parameters at the Radical Points}
\label{Local Parameters at the Radical Points}
\bigskip

For a radical ideal~$I$ having an $\mathcal{O}$-border basis,
we shall call the corresponding point of~$\mathbb{C}_{\mathcal{O}}$ a
{\bf radical point}. In the following, we want to construct 
explicit local parameters
for~$\mathbb{C}_{\mathcal{O}}$ near its radical points.
As a consequence, we shall recover the well-known facts 
that~$\mathbb{C}_{\mathcal{O}}$ is smooth of dimension~$\mu n$ 
at these points, and that it is a rational variety.

The idea of our construction is to use the complete intersection
representation of a radical ideal~$I$ having an $\mathcal{O}$-border
basis which is provided by the Shape Lemma (cf.~\cite{KR1}, Theorem~3.7.25).
Recall that a 0-dimensional ideal~$I$ is said to be in
{\bf normal $\ell$-position} for some $\ell\in P_1$ if we have
$\ell(p)\ne \ell(q)$ for distinct points $p,q\in\Supp(\mathcal{Z}(I))$.
Here the zero scheme $\mathcal{Z}(I)$ of~$I$ is defined over the
algebraic closure~$\overline{K}$ of~$K$.

\begin{proposition}
Let $I$ be a 0-dimensional radical ideal in~$P$ which has an
$\mathcal{O}$-border basis.
Assume that~$K$ has at least $\binom{\mu}{2}+1$ elements.

\begin{enumerate}
\item[(a)] It is possible to chose $\ell\in P_1$ such that~$\mathcal{Z}(I)$
is in normal $\ell$-position. 

\item[(b)] Write $\ell=\ell_1 x_1+\cdots +\ell_n x_n\in P_1$
with $\ell_1,\dots,\ell_n\in K$ and assume that $\ell_n\ne 0$.
Then we have $P/I\cong K[\ell]$ and the
minimal polynomial of~$\bar\ell$ in~$P/I$ is of the form $\chi(\ell)=
\ell^\mu-\lambda_{\mu}\ell^{\mu-1}-\cdots-\lambda_2\ell-\lambda_1$
where $\lambda_i\in K$. 

\item[(c)] The ideal~$I$ has a representation
$$
I=\bigl( \, \chi(\ell),\, x_1 - f_1(\ell),\, \cdots, x_{n-1} - f_{n-1}(\ell)
\,\bigr)
$$
where the polynomials $f_i(\ell)\in K[\ell]$ have degree $\le\mu-1$.

\end{enumerate}
\end{proposition}

\begin{proof} For claim~(a), see~\cite{KR1}, Proposition 3.7.22.
Claim~(b) follows from~\cite{KR1}, Theorem 3.7.23, and~(c) 
is the version of the Shape Lemma given in~\cite{KR1}, Theorem 3.7.25.
\end{proof}

Using the terminology of Section~2, the set $\mathcal{P}=\{1,\ell,\dots,
\ell^{\mu-1}\}$ is a pseudo order ideal because it is the image of the
order ideal $\{1,x_n,\dots,x_n^{\mu-1}\}$ under the linear change of
coordinates $\phi: P\longrightarrow P$ given by $\phi(x_i)=x_i$
for $i=1,\dots,n-1$ and $\phi(x_n)=\ell$. Next we define a grading
by $\deg_W(x_i)=\mu$ for $i=1,\dots,n-1$ and $\deg_W(x_n)=1$,
and we choose a term ordering~$\sigma$ which is compatible with this grading.
Then the set $\{1,x_n,\dots,x_n^{\mu-1}\}$ is a $\sigma$-cornercut and
its border is the set $\{x_1,\dots,x_{n-1},x_n^\mu\}$. Hence the set
$b\mathcal{P}=\{x_1,\dots,x_{n-1},\ell^\mu \}$ is the pseudo border
of~$\mathcal{P}$. 

Another way of stating the last claim of the proposition is to say that
the set $H=\{\chi(\ell), x_1-f_1(\ell),\dots,x_{n-1}-f_{n-1}(\ell)\}$
is a pseudo $\mathcal{P}$-border basis of~$I$. Pseudo border
bases of this shape are parametrized by $n\mu$ coefficients, namely
the coefficients $\lambda_1,\dots,\lambda_\mu$ of~$\chi$ and the
$(n-1)\mu$ coefficients of $f_1,\dots,f_{n-1}$. As $n\mu$ is the
dimension of~$\CO$ at the point corresponding to~$I$ (see~\cite{MS},
18.32), we shall now use the base change technique of Section~2 to parametrize
the principal component locally as follows. A similar result is shown
in~\cite{KLP} using a different technique.

\begin{theorem}\label{localparams}
Let $I$ be a 0-dimensional radical ideal in~$P$ which has an
$\mathcal{O}$-border basis. Suppose that there exist a linear
form $\ell=\ell_1 x_1 + \cdots +\ell_n x_n$ with $\ell_i\in K$
such that $\ell_n\ne 0$ and polynomials $\chi(\ell),f_i(\ell)\in K[\ell]$ 
such that
$$
I=\bigl( \, \chi(\ell),\, x_1 - f_1(\ell),\, \cdots, x_{n-1} - f_{n-1}(\ell)
\,\bigr)
$$
and such that
$\chi(\ell)=\ell^\mu -\lambda_\mu \ell^{\mu-1}-\cdots -\lambda_2\ell-\lambda_1$
and ${\deg(f_i(\ell))<\mu}$.
For every tuple $d=(d_{ij})\in K^{n\mu}$, we define
$\tilde\chi(\ell)=\ell^\mu - \sum_{i=1}^\mu(\lambda_i+d_{ni})\ell^{i-1}$
and $\tilde f_i(\ell)=f_i(\ell) +\sum_{j=1}^\mu d_{ij}\ell^{j-1}$
where $i\in\{1,\dots,n-1\}$. Then we let
$$
I_d= \bigl(\, \tilde\chi(\ell),\, x_1-\tilde f_1(\ell),\,\dots,\,
x_{n-1}-\tilde f_{n-1}(\ell) \,\bigr)
$$
\begin{enumerate}
\item[(a)] For all tuples~$d$ in a non-empty Zariski open neighborhood~$U$
of~$(0,\dots,0)$ in~$K^{n\mu}$, the ideal~$I_d$ has an $\mathcal{O}$-border basis
which can be computed by viewing $H_d=\{\tilde\chi(\ell), x_1-\tilde f_1(\ell),
\dots,x_{n-1}-\tilde f_{n-1}(\ell)\}$ as a pseudo $\mathcal{P}$-border basis
of~$I_d$ with respect to $\mathcal{P}=\{1,\ell,\dots,\ell^{\mu-1}\}$
and by applying Corollary~\ref{basechangeforI} to it.

\item[(b)] The morphism $\Gamma: U\longrightarrow \mathbb{C}_{\mathcal{O}}$
given by $d\mapsto I_d$ yields an isomorphism between the open set~$U$
in~$K^{n\mu}$ and the open set~$\mathbb{B}_{\mathcal{O},\mathcal{P}}\cap \CO$
in~$\CO$. In particular, the variety~$\mathbb{C}_{\mathcal{O}}$ is rational, and
it is smooth of dimension $n\mu$ at its radical points.
\end{enumerate}
\end{theorem}

\begin{proof}
Let us consider the pseudo order ideal $\mathcal{P}=\{1,\ell,\dots,\ell^{\mu-1}\}$
as the image of the cornercut $\{1,x_1,\dots,x_n^{\mu-1}\}$ under a
linear change of coordinates.
Then its pseudo border is $b\mathcal{P}=\{x_1,\dots,x_{n-1},\ell^\mu\}$,
and $H=\{\chi(\ell), x_1-f_1(\ell),\dots,x_{n-1}-f_{n-1}(\ell)\}$
is the pseudo $\mathcal{P}$-border basis of~$I$.

First we prove~(a). The ideal~$I$ has both an $\mathcal{O}$-border basis 
and a pseudo $\mathcal{P}$-border basis. For every $d\in K^{n\mu}$, the
set~$H_d$ is a pseudo $\mathcal{P}$-border basis of~$I_d$.
By Theorem~\ref{basechange}.a, there is
a non-empty Zariski open neighbourhood~$U$ of~$(0,\dots,0)$ in~$K^{n\mu}$
such that $I_d$ has an $\mathcal{O}$-border basis for all $d\in U$. 
It is the open set defined by~$\det(M_D)\ne 0$ where
$M_D \in\Mat_\mu(K[d_{ij}])$ is the matrix
such that $\mathcal{O}\equiv \mathcal{P}\cdot M_D\; (\mod\, I_d)$
where we view the elements $d_{ij}$ as indeterminates.
For all $d\in U$, we can use Corollary~\ref{basechangeforI} to compute the
$\mathcal{O}$-border basis of~$I_d$. For this purpose, we have to calculate
matrices $M_d\in\Mat_\mu(K)$ and $N_d\in \Mat_{\mu,\nu}(K)$
such that $\mathcal{O}=\mathcal{P}\cdot M_d$ and
$\partial\mathcal{O}=\mathcal{P}\cdot N_d$ hold in~$P/I_d$.
Then the matrix representing the ideal~$I_d$ in~$\BO$
is given by $C_{I_d}=(M_d)^{-1}\cdot N_d$.

It remains to prove~(b). The morphism~$\Gamma$ sends a tuple $d=(d_{ij})$
to the point represented by the matrix $C_{I_d}=(M_d)^{-1}\cdot N_d$. 
We claim that it maps the open set $U=K^{n\mu}\setminus \mathcal{Z}(\det(M_C))$
isomorphically to the open set $V={\mathbb{B}_{\mathcal{O},\mathcal{P}}\cap \CO}$
in~$\CO$. Given a matrix~$C_J$ representing a point in~$V$,
we know that the corresponding ideal~$J$ has both an $\mathcal{O}$-border
basis and a pseudo $\mathcal{P}$-border basis. Using Corollary~\ref{basechangeforI}
again, we can compute a tuple $d\in U$ such that $J=I_d$, and this
yields a morphism from~$V$ to~$U$ which inverts~$\Gamma$.
\end{proof}

It is well-known that~$\BO=\CO$ is non-singular for the case of
$n=2$ indeterminates and that the following example shows that~$\CO$
is not always non-singular at its non-radical points.

\begin{example}
Let $P=\mathbb{Q}[x,y,z]$, and let $\mathcal{O}=\{1,x,y,z\}$.
Then $\BO=\CO$ is a scheme of dimension~12 in~$\mathbb{A}^{24}$.
We compute $I(\BO)$ and see that we can project~$\BO$
isomorphically to an 18-dimensional affine space by
eliminating $c_{11},\dots,c_{16}$ (cf.~\cite{KR3}). The result
is a variety $\pi(\BO)\subset \mathbb{A}^{18}$
whose vanishing ideal is generated by~15 homogeneous
polynomials of degree~2.
Hence its vertex $(0,\dots,0)$ is singular. The corresponding ideal
is the border term ideal
$\langle\partial\mathcal{O}\rangle=\langle x^2,xy,xz,y^2,yz,z^2\rangle$.
\end{example}

In~\cite{Ha2}, Sect.~2, it is explained that~$\CO$ can be realized as
the blowup of the Chow variety ${\rm Spec}(K[\y1,\dots,\ymu]^{S_\mu}$
at an explicitly given ideal. A different construction for~$\CO$, permitting
similar conclusions, is contained in~\cite{ES}. Although it is in principle
possible from these constructions to obtain local parameters for~$\CO$
at its radical points, we believe that our construction is
more elementary and explicit.

Let us compare the results of Section~2 to the
preceding theorem. Given an arbitrary $K$-rational point of~$\BO$,
i.e.\ an ideal~$I$ having an $\mathcal{O}$-border basis, we can
use linear changes of coordinates as in Proposition~\ref{linchange}
to construct a flat family of ideals whose base space is an open subset of an
$n(n+1)$-dimensional affine space and whose special fiber is~$I$.
However, this is in general much smaller than the local dimension of~$\BO$ at
the point representing~$I$. If~$I$ is reduced, Theorem~\ref{localparams}
allows us to do much better: we construct a flat family over a $n\mu$-dimensional
base space, and this is precisely the local dimension of~$\BO$
at the point representing~$I$.

An application of the preceding theorem is the possibility
to connect two arbitrary radical ideals having $\mathcal{O}$-border
bases via a sequence of two explicit flat families. 

\begin{remark}\label{connect2rad}
Let $K$ be an infinite field, let $P=K[x_1,\dots,x_n]$,
let $I,I'\subset P$ be two radical ideals which have
$\mathcal{O}$-border bases, and let $c_I,c_{I'}$ be the points
in~$\CO$ representing them.

For a generically chosen $\ell\in P_1$, the hypotheses of
the theorem are satisfied with respect to~$I$ and~$I'$
(see \cite{KR1}, Prop.\ 3.7.22).
By part~(c) of the theorem, there exist an open
neighborhood~$U$ of~$c_I$ in~$\CO$ and an open
neighborhood~$U'$ of~$c_{I'}$ in~$\CO$ such that
the restriction of the universal flat family to~$U$ and~$U'$
is given by explicitly defined morphisms.

Since~$\CO$ is irreducible, there exists a point $c_J\in U\cap
U'$ representing a radical ideal~$J$. Both~$U$ and~$U'$
are isomorphic to open subsets of~$\mathbb{A}^{n\mu}$. 
Our task is to find an explicit flat family
connecting~$c_I$ and~$c_J$. In an analogous way, we can then
connect~$c_{I'}$ and~$c_J$.

The points $p_I=\Gamma^{-1}(c_I)$ and $p_J=\Gamma^{-1}(c_J)$
are contained in an open subset~$U$ of the affine space $\mathbb{A}^{n\mu}$.
By Theorem~\ref{basechange}, we have an explicit polynomial~$F$ 
whose non-vanishing defines this open set. Thus we can connect
the points~$p_I$ and~$p_J$ by a line~$L$ and get a
$K$-algebra homomorphism $K[d_{ij}]_F \longrightarrow K[z]_f$
which represents the inclusion $(L\cap U)\subseteq U$. Here
$f\in K[z]\setminus\{0\}$ defines the points in $L\setminus U$.

After applying~$\Gamma$, we get an explicit punctured rational curve
$\Psi: C_{\mathcal{O}}\longrightarrow K[z]_f$ which
connects~$c_I$ to~$c_J$ in~$\mathbb{B}_{\mathcal{O},\mathcal{P}}\cap \CO$. 
Then, by restricting the universal flat family
$\Phi: B_{\mathcal{O}} \longrightarrow U_{\mathcal{O}}$
to this punctured rational curve, we find a
flat family deforming~$P/I$ to~$P/J$.
\end{remark}

Another application of Theorem~\ref{localparams} is the
possibility to connect an arbitrary radical point of~$\BO$
to the {\bf monomial point} $o=(0,\dots,0)$ representing the monomial
ideal $\langle b_1,\dots,b_\nu\rangle$ via explicitly
defined flat families. For this application we need one further
ingredient, namely distractions, which we are now going to
recall from~\cite{KR2}, Section~6.2.
To simplify the discussion, let us assume that the field~$K$
has sufficiently many elements.

\begin{definition}
For $i=1,\dots,n$, let $\pi_i=(c_{i1},c_{i2},\dots)$
be tuples consisting of sufficiently many pairwise distinct elements
of~$K$.

\begin{enumerate}
\item For a term $t=x_1^{\alpha_1}\cdots x_n^{\alpha_n}$, the
polynomial $D_\pi(t)=\prod_{i=1}^n \prod_{j=1}^{\alpha_i}
(x_i-c_{ij})$ is called the {\bf distraction} of~$t$ with respect
to $\pi=(\pi_1,\dots,\pi_n)$.

\item Let $I=\langle b_1,\dots,b_\nu\rangle$ be the border term ideal
of~$\mathcal{O}$, and let $c\mathcal{O}= \{ c_1,\dots,c_r\}$ be the
corner set of~$\mathcal{O}$, i.e.\ the set of minimal monomial generators of~$I$.
Then the ideal $D_\pi(I)=
\langle D_\pi(c_1),\dots,D_\pi(c_r)\rangle$ is called the
{\bf distraction} of~$I$ with respect to~$\pi$.
\end{enumerate}
\end{definition}

The following properties of the distraction of the border term ideal
are shown in~\cite{KR2}, Theorem 6.2.12.

\begin{proposition}\label{distract}
Let $\pi=(\pi_1,\dots,\pi_n)$ be chosen as in the preceding definition,
and let $c\mathcal{O}=\{c_1,\dots,c_r\}$ be the corner set
of the border term ideal~$I$ of~$\mathcal{O}$.

\begin{enumerate}
\item The distraction $D_\pi(I)$ is a radical ideal.

\item For every term ordering~$\sigma$, the set
$\{D_\pi(c_1),\dots,D_\pi(c_r)\}$ is the reduced
$\sigma$-Gr\"obner basis of~$D_\pi(I)$.
In particular, we have $\LT_\sigma(D_\pi(I))=I$
and $\mathbb{T}^n\setminus \LT_\sigma(I)=\mathcal{O}$.

\item The ideal $D_\pi(I)$ is 0-dimensional and has
an $\mathcal{O}$-border basis.
\end{enumerate}
\end{proposition}

Now we are ready to connect any radical point of~$\BO$
to the monomial point via three explicit flat families.

\begin{remark}\label{connect2mon}
Let $K$ be an infinite field.
Suppose we are given a reduced 0-dimensional ideal~$I$
in~$P$ which has an $\mathcal{O}$-border bases.
To find explicit flat families connecting $P/I$ to 
$P/\langle b_1,\dots,b_\nu\rangle$, we proceed as follows.

\begin{enumerate}
\item For $i=1,\dots,n$, choose tuples~$\pi_i$ of sufficiently many
distinct elements of~$K$. For every $\lambda\in K$, let $\lambda\pi =
(\lambda\pi_1,\dots,\lambda\pi_n)$. Then form the family
$\Pi: \mathbb{A}^1 \longrightarrow \CO$ defined by $\lambda\mapsto
D_{\lambda\pi}(\langle b_1,\dots,b_\nu\rangle)$.
By Proposition~\ref{distract}, there exists a family of polynomials
which yields a Gr\"obner basis of each fiber of this family. Hence
$\Pi$ is a flat deformation connecting the border term ideal to the
radical ideal $J=D_\pi(\langle b_1,\dots,b_\nu\rangle)$.

\item Now use Remark~\ref{connect2rad} to find two explicit 
flat families connecting $P/J$ to~$P/I$.
\end{enumerate}
\end{remark}

Notice that by using the method of the previous remark we may not
always find a non-punctured rational curve in~$\CO$
connecting $c_I$ to the monomial point~$o$, although such a curve
might exist.

We end this section with some examples which illustrate
the construction of the explicit flat families in Remarks~\ref{connect2rad}
and~\ref{connect2mon}. In the first one we connect the points $c_I$ and
$c_J$ corresponding to two radical ideals by a rational curve in~$\BO$,
but one point of the rational curve is not contained in the open
set $\Gamma(U)=\mathbb{B}_{\mathcal{O},\mathcal{P}}\cap \CO$
of Theorem~\ref{localparams}.b. 

\begin{example}
Let $K=\mathbb{Q}$, let $P=\mathbb{Q}[x,y]$, and let $\mathcal{O}=
\{1,x,y,xy\}$. The vanishing ideals of the two point sets
$\mathbb{X}=\{(0,0),\, (0,1),\, (1,-1),\, (1,2)\}$ and
$\mathbb{Y}=\{(0,0),\, (0,1),\, (1,-1),\, (-1,2)\}$ both have
$\mathcal{O}$-border bases, namely
\begin{eqnarray*}
\mathcal{I}(\mathbb{X}) &=& ( y^2-2x-y,\, x^2-x,\, xy^2-2x-xy,\, x^2y-xy )\\
\mathcal{I}(\mathbb{Y}) &=& ( y^2-\tfrac{2}{3}x+y+\tfrac{4}{3}xy,\ x^2-\tfrac{1}{3}x
+\tfrac{2}{3}xy,\, xy^2-2x-xy,\, x^2y+\tfrac{4}{3}x+\tfrac{1}{3}xy )
\end{eqnarray*}
Now we use the explicit description of~$\BO$ worked out in~\cite{KR3}, Example~3.8.
It provides an isomorphism $\Gamma: \mathbb{A}^8 \longrightarrow \BO$
which corresponds to the $\mathbb{Q}$-algebra homomorphism
$$ 
B_{\mathcal{O}} \;\longrightarrow\; \mathbb{Q}[c_{21},c_{23},c_{32},c_{34},
c_{41},c_{42},c_{43},c_{44}] 
$$
given by $(c_{11},c_{12},\dots,c_{44}) \mapsto (c_{21},c_{23},\dots,c_{44})$.
Under this isomorphism, the point representing $\mathcal{I}(\mathbb{X})$
corresponds to $(c_{21},\dots,c_{44})=(2,2,0,0,0,0,1,1)=p_1$, and the point 
representing $\mathcal{I}(\mathbb{Y})$ corresponds to $p_2=(\tfrac{2}{3},2,0,0,
-\tfrac{4}{3},-\tfrac{2}{3},1,-\tfrac{1}{3})$.

To connect~$p_1$ and~$p_2$, we use the line $L=\{(-1+2a)p_1+
(\tfrac{1}{2}-\tfrac{1}{2}a)p_2 | a\in \mathbb{Q} \}$. In this way we 
get the point $p_1$ for $a=1$ and the point $p_2$ for $a=-1$.
For the corresponding family of ideals, we use
$$
c_{21}=\tfrac{2}{3}a+\tfrac{4}{3}, c_{23}=2, c_{32}=c_{34}=0,
c_{41}=\tfrac{2}{3}a-\tfrac{2}{3}, c_{42}=\tfrac{1}{3}a-\tfrac{1}{3},
c_{43}=1, c_{44}=\tfrac{2}{3}a+\tfrac{1}{3}
$$
This yields the family of ideals whose border bases are represented
parametrically by
\begin{eqnarray*}
G_a &=& \{ y^2 - (-\tfrac{2}{3}a^2+\tfrac{2}{3})- (\tfrac{2}{3}a+\tfrac{4}{3})x
- (-\tfrac{2}{3}a^2+\tfrac{5}{3})y - (\tfrac{2}{3}a-\tfrac{2}{3})xy,\\
&& x^2 - (\tfrac{1}{3}a +\tfrac{2}{3})x -(\tfrac{1}{3}a-\tfrac{1}{3})xy,\\
&& xy^2-2x-xy,\\
&& x^2y - (\tfrac{2}{3}a-\tfrac{2}{3})x - (\tfrac{2}{3}a+\tfrac{1}{3})xy \}
\end{eqnarray*}
The ideal $I_a$ generated by~$G_a$ satisfies 
$$
I_0 = (x-1,\, y+1) \;\cap\; (x,\,y+\tfrac{1}{3}) \;\cap\; (x^2,\,y-2)
$$
in the case $a=0$ and
$$
I_a = (x-1,\,y+1) \;\cap\; (x-a,\, y-2) \;\cap\; (x,\, \tfrac{3}{2}y^2
+(a^2-\tfrac{5}{2})y + (a^2-1))
$$
for $a\ne 0$. Thus the ideal~$I_0$ is not reduced, 
but all other ideals of the family are.

Geometrically, this can be explained as follows. No set of points in the 
family can have three points on the line $x=0$, since then the polynomial 
$xy+x$ vanishes on all four points and there is no $\mathcal{O}$-border basis.
Therefore the point $(0,1)$ ``moves up'' and helps the point $(a,2)$ to get 
across this line by forming a non-reduced scheme. 
\end{example}

The punctured rational curves constructed in Remark~\ref{connect2rad} 
may sometimes be restrictions
of (non-punctured) rational curves on the Hilbert scheme
whose special points lie outside the border basis scheme,
not just outside the set~$\Gamma(U)$. A case in point is Example~3.9
of~\cite{KR3} where the value $a=0$ corresponds to an ideal~$I_0$
which has no $\mathcal{O}$-border basis.

Our last example shows that the flat families we construct
do, in general, not preserve the geometry of the corresponding
sets of points. 

\begin{example}
In the setting of the preceding example, we replace the scheme~$\mathbb{Y}$ 
by $\mathbb{Y}'=\{(0,0),\, (-1,1),\, (1,-1),\, (1,2)\}$. 
Arguing as above, we see that~$\mathbb{Y}'$ is represented in~$\mathbb{A}^8$
by the point $p_3=(2,2,1,1,-1,-1,0,0)$. Now we connect~$p_1$ and~$p_3$ by a line
such that~$p_1$ corresponds to the parameter value $a=0$ and~$p_3$ corresponds to
$a=1$. The resulting family of ideals uses 
$$
c_{21}= c_{23}=2,\; c_{32}=c_{34}=a,\; c_{41}=c_{42}=-a,\; c_{43}=c_{44}=1-a
$$
and its border bases are parametrically given by
\begin{eqnarray*}
H_a &=& \{ y^2  - 2x +(a^2-2a-1)y +axy, \\
&& x^2 - x - ay +a xy, \\
&& xy^2 - 2x +(a^2-2a)y -(1-a) xy, \\
&& x^2y - ay - (1-a)xy \}
\end{eqnarray*}

The ideal $J_a$ generated by $H_a$ is reduced for every~$a$.
For $a=0$, the ideal~$J_0=\mathcal{I}(\mathbb{X})$ corresponds to
a set of four points in general position, but for $a=1$ the ideal
$J_1=\mathcal{I}(\mathbb{Y}')$ corresponds to four points, three of which
lie on the line $\mathcal{Z}(x+y)$. In geometrical jargon one can express this 
by stating that the set~$\mathbb{X}$ has the Cayley-Bacharach property, 
but~$\mathbb{Y}'$ doesn't. In this sense the flat family
did not preserve the geometry of the point set.
\end{example}

\bigbreak
%%%%%%%%%%%%%%%%%%%%%%%%%%%%%%%%%%%%%%
\subsection*{Acknowledgements}
The first author thanks the Dipartimento di Matematica
of Universit\`a di Genova for the warm hospitality he enjoyed
during part of the preparation of this paper. Likewise, the second
author is grateful to the Fakult\"at f\"ur Informatik und Mathematik
of Universit\"at Passau for the hospitality he enjoyed during
other parts of the writing of this paper.
The examples in this paper were computed using the computer
algebra system \cocoa\ (see~\cite{CoCoA}) and the border
basis package of \apcocoa\ (see~\cite{ApCoCoA}).

\bigbreak
%%%%%%%%%%%%%%%%%%%%%%%%%%%%%%%%%%%%%%
%   Bibliography
%%%%%%%%%%%%%%%%%%%%%%%%%%%%%%%%%%%%%%

\end{document}